\documentclass[12pt]{article}
\usepackage{amsmath,amssymb,amsfonts,latexsym}

\newtheorem{theorem}{Theorem}
\newtheorem{lemma}{Lemma}

\newtheorem{proposition}{Proposition}
\newtheorem{conjecture}{Conjecture}

\newtheorem{example}{Example}

\begin{document}

\title{On the Expected Value of the Minimum Assignment}
\author{Marshall W. Buck \and Clara S. Chan \and David P. Robbins
\thanks{All three authors are affiliated with the Center for Communications Research, Princeton, NJ 08540.  Email: buck@idaccr.org, clara@idaccr.org, robbins@idaccr.org}}
\date{May 10, 2000}
\maketitle

\begin{abstract}
The minimum $k$-assignment of an $m\times n$ matrix $X$ is the minimum sum
of $k$ entries of $X$, no two of which belong to the same row or column.
If $X$ is generated by choosing each entry independently from the exponential
distribution with mean 1, then Coppersmith and Sorkin conjectured that the
expected value of its minimum $k$-assignment is
$$\sum_{i,j \ge 0, \;i+j<k} \frac{1}{(m-i)(n-j)}$$
and they (with Alm) have proven this for $k\le 4$ and in certain cases when
$k=5$ or $k=6$.
They were motivated by the special case of $k=m=n$, where the expected value
was conjectured by Parisi to be
$\sum_{i=1}^k \frac{1}{i^2}.$
In this paper we describe our efforts to prove the Coppersmith--Sorkin
conjecture.  We give evidence for the following stronger conjecture, which
generalizes theirs.

\textbf{Conjecture} Suppose that $r_1,\dots,r_m$ and $c_1,\dots,c_n$ are positive real numbers.
Let $X$ be a random $m \times n$ matrix in which entry $x_{ij}$ 
is chosen independently from the exponential distribution with mean $\frac1{r_ic_j}$.
Then the expected value of the minimum $k$-assignment of $X$ is 
$$\sum_{I,J}
  (-1)^{ k - 1 - |I| - |J| } \cdot 
   \binom{m + n - 1 - |I| - |J|}{k - 1 - |I| - |J|}
      \frac{1}{( \sum_{i \notin I}r_i) \cdot ( \sum_{j \notin J} c_j )}.$$
Here the sum is over proper subsets $I$ of $\{1,\dots,m\}$ and
$J$ of $\{1,\dots,n\}$ whose cardinalities $|I|$ and $|J|$ satisfy 
$|I|+|J|<k$.
\end{abstract}

\section{Problem Description and Background}\label{sec:intro}

Suppose that $k$, $m$ and $n$ are positive integers with $k \le m\le n$.
A {\em minimum $k$-assignment\/}\ of an $m \times n$ matrix $X$ is a set of 
$k$ entries of $X$, no two of which belong to the same row or column,
whose sum is as small as possible.  We denote the value of this minimum
sum by ${\min}_k(X)$.  

We say that a random real number $x$ is
{\em exponentially distributed\/}\ with {\em rate\/}\ $a$ if
it is chosen according to the density $ ae^{-ax}$, $x\ge 0$.  The mean
value of a rate $a$ exponentially distributed quantity is $1/a$.

Suppose that we
generate a random $m \times n$ matrix $X$ by
choosing each entry independently from the exponential
distribution with rate 1. In \cite{cs} Coppersmith and Sorkin conjectured that the
expected value of its minimum $k$-assignment is
\begin{conjecture}\label{conj:cs}
\begin{equation}\label{eq:cs}
E({\min}_k(X))=\sum_{i,j \ge 0, \;i+j<k} \frac{1}{(m-i)(n-j)}.
\end{equation}
\end{conjecture}
In \cite{as} Alm and Sorkin show that this conjecture is correct when
$k\leq4$, when $k=m=5$, and when $k=m=n=6$.

The conjecture of Coppersmith and Sorkin generalized a conjecture of 
Parisi \cite{p} who considered the case
$k=m=n$.  In this case, as shown in \cite{cs},  (\ref{eq:cs})  reduces 
to 
\begin{equation}\label{eq:par}
E({\min}_k(X))=\sum_{i=1}^k \frac{1}{i^2}.
\end{equation}

In this paper we describe our efforts to prove these conjectures.
Our main result is Conjecture \ref{conj:r1}, which
generalizes the Coppersmith-Sorkin conjecture.

We will say that a matrix $X$ is {\em random exponential\/}\  with
{\em rate matrix\/}\ $A=(a_{ij})$ if each entry $x_{ij}$ is chosen
independently according to the exponential
distribution with rate $a_{ij}$.  The expected value
of the minimum $k$-assignment of such a matrix $X$ is then
a function of the rate matrix $A$.  We denote this function by $E_k(A)$. 

We will show that $E_k(A)$ is a rational function of the
rates $a_{ij}$ and give an explicit method for computing it,
at least in principle.  Then we will specialize to the case
when the rate matrix has rank 1, for which we have the
following explicit formula.
\begin{conjecture}\label{conj:r1}
Suppose that $r_1,\dots,r_m$ and $c_1,\dots,c_n$ are positive real
numbers and that $a_{ij}=r_ic_j$.  Let $X$ be a random $m \times n$
matrix in which entry $x_{ij}$ is chosen independently from the
exponential distribution with rate $a_{ij}$.  Then the expected value
of the minimum $k$-assignment of $X$ is
$$\sum_{I,J} (-1)^{ k - 1 - |I| - |J| } \cdot
  \binom{m + n - 1 - |I| - |J|}{k - 1 - |I| - |J|}
  \frac{1}{( \sum_{i \notin I}r_i) \cdot ( \sum_{j \notin J}c_j )}.$$
Here the sum is over proper subsets $I$ of $\{1,\dots,m\}$
and $J$ of $\{1,\dots,n\}$ whose cardinalities $|I|$ and $|J|$ satisfy
$|I|+|J|<k$.
\end{conjecture}

\begin{example}\label{ex:1}
The expected value of the minimum 1-assignment of a random exponential
matrix with rate matrix $a_{ij}=r_ic_j$ is 
\[
\frac{1}{( \sum_{i}r_i) \cdot ( \sum_{j} c_j )}.
\]
\end{example}
\begin{example}
The expected value of the minimum 2-assignment of a $3\times3$
random exponential matrix with rate matrix $a_{ij}=r_ic_j$ is
\begin{eqnarray*}
&&\left(\frac1{r_2+r_3}+\frac1{r_1+r_3}+\frac1{r_1+r_2}\right)
\frac1{c_1+c_2+c_3} \\
&+&\left(\frac1{c_2+c_3}+\frac1{c_1+c_3}+\frac1{c_1+c_2}\right)
\frac1{r_1+r_2+r_3} \\
&-&\frac5{(r_1+r_2+r_3)(c_1+c_2+c_3)}
\end{eqnarray*}
\end{example}

We will provide evidence in support of Conjecture~\ref{conj:r1}.

We also have a stronger conjecture for which we will provide evidence,
although perhaps this evidence is not as strong as that for
Conjecture~\ref{conj:r1}.
A matrix can have several minimum $k$-assignments for some value of $k$.
However, with probability 1, a random matrix has a single 
minimum $k$-assignment for each $k$.  Suppose that $M$ is a
$(k-1) \times (k-1)$ submatrix
of $X$ and that $\chi_M(X)$ is the function with value 1 when $M$
contains a minimum $(k-1)$-assignment of $X$ and 0 otherwise.  Then we define
the {\em expected contribution\/}\ of $M$ to the minimum $k$-assignment of $X$
as the expected value of the random variable $\chi_M(X) {\min}_k(X)$.  It is
clear that $E({\min}_k(X))$ is the sum of the expected contributions of all the
$(k-1) \times (k-1)$ submatrices.  Our stronger conjecture gives a formula
for the expected contribution of $M$ when the rate matrix has rank 1.

\begin{conjecture}\label{conj:contrib}
Suppose that $A=(r_ic_j)$ is a positive $m \times n$ matrix with rank 1
and that $X$ is a random exponential matrix with rate matrix $A$.  Let
$I$ be a set of $k-1$ elements of $\{1,\ldots,m\}$, let $J$ be a set of
$k-1$ elements of $\{1,\ldots,n\}$, and let $M$ be the $(k-1)\times (k-1)$
submatrix of $X$ with rows indexed by $I$ and columns indexed by $J$.
Then the expected value of $\chi_M(X) \cdot {\min}_k(X)$ is
$$
\sum_{i,j}
\left(\prod_{t=1}^{k-1} \frac{ r_{i_t}c_{j_t}}
 {\left(R-\sum_{s=1}^{t-1} r_{i_s}\right)\left(C-\sum_{s=1}^{t-1} c_{j_s}\right)}\right)
\sum_{t,u\ge 0, t+u < k}
\frac{1}{(R-\sum_{s=1}^tr_{i_s})(C-\sum_{s=1}^u c_{j_s})}
$$
where the outer sum is over permutations $(i_1,\dots,i_{k-1})$ and
$(j_1,\dots,j_{k-1})$ of $I$ and $J$, respectively, 
and $R$ and $C$ denote the sums of all $r_i$'s and all $c_j$'s, respectively.

\end{conjecture}

We shall see that Conjecture~\ref{conj:contrib} implies Conjecture~\ref{conj:r1} and
that Conjecture~\ref{conj:r1} in turn implies Conjecture~\ref{conj:cs}.

Section~\ref{sec:general} discusses what we know for the
expected minimum assignment when the rate matrix is arbitrary.

In Section~\ref{sec:rank1} we discuss the way we arrived at Conjecture~\ref{conj:r1}
and give some equivalent formulations, one of which is directly implied by 
Conjecture~\ref{conj:contrib}. 

We discuss the computational evidence for our conjectures in Section~\ref{sec:comp}.

Section~\ref{sec:extra} gives additional evidence for Conjecture~\ref{conj:r1}.

We would like to thank Jim Propp for bringing this problem to our attention.

\section{Theory for a general rate matrix}\label{sec:general}

\subsection{Expected value for a general rate matrix}

We begin by showing that the general formula for the
expected value of the minimum assignment of a random exponential
matrix is a rational function of the rates, with denominators
factoring into linear terms of special form.

Recall that $k,m,n$ are positive integers with $k\le m \le n$ and that
$A=(a_{ij})$ is a positive $m \times n$ matrix.  We form a random 
matrix $X$ by choosing $x_{ij}$ independently from the exponential
distribution with rate $a_{ij}$.  The expected value of the minimum
$k$-assignment of $X$ is then a function of $A$, which we will denote
by $E_k(A)$.

By definition of expected value, $E_k(A)$ is given by the integral expression 
$$E_k(A)=\left(\prod_{i,j}a_{ij}\right) \int_{X \ge 0} {\min}_k(X) e^{-A\cdot X} dX$$
where the integral is taken over the space of all nonnegative
matrices $X$.  Here $A\cdot X$ denotes the dot product
$\sum_{i,j} a_{ij}x_{ij}$ and $dX$ denotes the product $\prod_{i,j} dx_{ij}$.

We denote by ${\cal S}_k$ the set of all $m \times n$ matrices $\sigma$
such that all the entries of $\sigma$ are 0's except for $k$
entries which are 1's, no two in the same row or column.  There are
$k!\binom{m}{k}\binom{n}{k} $ such matrices in ${\cal S}_k$ 
and these we identify in
the obvious way with the possible locations of the minimum
$k$-assignment of $X$.  In particular,
$$  {\min}_k(X)=\min_{\sigma \in {\cal S}_k} \left(\sigma \cdot X\right) \,. $$

For each $\sigma$ we denote by $P_\sigma$ the set of 
nonnegative matrices $X$ for which
${\min}_k(X)=\sigma \cdot X$; that is, $P_\sigma$ is the set of 
nonnegative matrices $X$ for which the
minimum $k$-assignment is $\sigma$.  Thus, we have
\begin{equation}\label{eq: integral0}
E_k(A)=\left(\prod_{i,j}a_{ij}\right)\sum_{\sigma\in {\cal S}_k}\left[ \int_{X \in P_\sigma} 
   (\sigma \cdot X) e^{-A\cdot X} dX\right] \,.
\end{equation}

Note that  each of the sets
$P_\sigma$ is a polyhedral cone determined by a finite set of homogeneous 
linear inequalities
$ \sigma \cdot X \le \tau \cdot X$ for all $\tau \in {\cal S}_k$.
As a consequence, each 
$P_\sigma$ can be decomposed into a finite collection ${\cal C}_\sigma$ of 
simplicial cones.  It seems difficult to give an explicit description of
${\cal C}_\sigma$.  Nevertheless, we can derive some useful properties of
$E_k(A)$ from the fact that this decomposition exists.
First we rewrite (\ref{eq: integral0}) as
$$E_k(A)=\left(\prod_{i,j}a_{ij}\right)\sum_{\sigma\in {\cal S}_k}\sum_{C \in {\cal C}_\sigma}
\left[ \int_C 
   (\sigma \cdot X) e^{-A\cdot X} dX\right] \,. $$

Each cone $C$ is
the set of nonnegative linear combinations of a set of $mn$ linearly independent vectors $V_i$, $i=1,\dots,mn$, where each $V_i$ is a nonnegative $m \times n$ matrix.
For the part of the integral over $C$, we make the substitution 
$X=\sum_i u_i V_i $, where $U=(u_1,\dots,u_{mn})$ ranges over all nonnegative
$mn$-tuples.  We can then explicitly compute the integral over $C$ as
\begin{eqnarray*}
\int_{C} (\sigma \cdot X) e^{-A\cdot X} dX & = &
|\det V| \int_{U \ge 0}\sum_{i=1}^{mn} u_i(\sigma \cdot V_i)
  e^{-\sum_{j=1}^{mn} u_j A\cdot V_j} dU \\
& = & 
|\det V| \sum_{i=1}^{mn}
\left[ (\sigma \cdot V_i)\int_{U \ge 0} u_i e^{-\sum_{j=1}^{mn} u_j A\cdot V_j} dU\right]\\
& = & 
| \det V| 
\left(\sum_{i=1}^{mn} \frac{\sigma \cdot V_i}{A\cdot V_i}\right)
\left(\prod_{i=1}^{mn} \frac{1}{A\cdot V_i}\right)
 \\
\end{eqnarray*}
where $|\det V|$ is the $mn$-volume of the parallelepiped 
determined by $V_1,\dots,V_{mn}$.
Thus, we obtain the expression
\begin{equation}\label{eq:2}
E_k(A)=\left(\prod_{i,j}a_{ij}\right)\sum_{\sigma\in {\cal S}_k}\sum_{C \in {\cal C}_\sigma}
| \det V| 
\left( \sum_{i=1}^{mn} \frac{\sigma \cdot V_i}{A\cdot V_i}\right) 
\left(\prod_{i=1}^{mn} \frac{1}{A\cdot V_i}\right)
\,. 
\end{equation}
Note that although the vectors $V_i$ depend on $C$ and $\sigma$, they do
not depend on $A$.  Thus, $E_k(A)$ is a rational function of the
$a_{ij}$'s,  homogeneous of degree $-1$.

We can obtain more information about the rational function $E_k(A)$ 
by constructing, for each $\sigma \in {\cal S}_k$,
a finite set of generators of $P_\sigma$ in the sense that every element
of $P_\sigma$ is a nonnegative linear combination of the
generators.

For this purpose we define two classes of matrices.
First, for any $1\le i\le m$ and $1\le j\le n$ we define
$e_{ij}$ to be the matrix
that is all zero except for a single 1 at position $(i,j)$.
Next, for any sets $I\subseteq\{1,\dots,m\}$ and $J\subseteq \{1,\dots,n\}$,
we define $V_{IJ}$ to be the matrix obtained from the all 1's matrix
by zeroing out all entries in the rows indexed by $I$ and the columns
indexed by $J$. 
It is easy to see that 
${\min}_k(V_{IJ})=\max(0,k-|I|-|J|)$.  Thus $V_{IJ}$ is in $P_\sigma$ if 
and only if $\sigma \cdot V_{IJ}=\max(0,k-|I|-|J|)$.

\begin{theorem}\label{th:1}
Every element of $P_\sigma$ is a nonnegative linear combination
of $e_{ij}$'s with $e_{ij} \cdot \sigma=0$ and
$V_{IJ}$'s in $P_{\sigma}$ with $|I|+|J|<k$.
\end{theorem}

We prove Theorem~\ref{th:1} using a reduction procedure on the matrices
of $P_\sigma$.  Let $X$ be a matrix in
$P_\sigma$ and suppose that ${\min}_k(X)=s$.
We choose an arbitrary linear order for the $e_{ij}$'s
and denote this ordered set by $e_1,e_2,\dots,e_{mn}$.
Then we choose a sequence of nonnegative real numbers
$\alpha_1,\alpha_2,\dots,\alpha_{mn}$ as follows.
Once $\alpha_1,\dots,\alpha_{i-1}$ are chosen we
select $\alpha_i$ as large as possible so that 
$X-(\alpha_1 e_1 + \alpha_2 e_2 +\cdots+ \alpha_i e_i)$ is nonnegative
and has minimum $k$-assignment with value $s$.  

Set 
$$Y=X-(\alpha_1 e_1 + \alpha_2 e_2 +\cdots+ \alpha_{mn} e_{mn}). $$

Note that if $e_i \cdot \sigma \ne 0$, then we
will have $\alpha_i=0$, since otherwise 
$\sigma \cdot(X-(\alpha_1 e_1 + \alpha_2 e_2 +\cdots+ \alpha_i e_i))<s.$
Thus $X$ is $Y$ plus a nonnegative linear combination of the $e_{ij}$'s given
in Theorem~\ref{th:1}.

We say that an entry $y_{ij}$ of a matrix $Y$ \emph{participates} in a minimum 
$k$-assignment if there is a minimum $k$-assignment using the entry $y_{ij}$.

We say that a nonnegative matrix $Y=(y_{ij})$ is \emph{$k$-reduced} if every nonzero
entry of $Y$ participates in a minimum $k$-assignment. 

It is straightforward to see that the matrix $Y$ resulting from our
reduction process applied to $X \in P_\sigma$ is $k$-reduced and that
$Y\in P_\sigma$.  It remains to show that 
every $k$-reduced matrix $Y$ with minimum $k$-assignment $\sigma$ is
a nonnegative linear combination of the appropriate $V_{IJ}$.  This will
require a series of preliminary results.

First we need a simple combinatorial lemma.

\begin{lemma}\label{lem:0}
Suppose that $T$ is a matrix all of whose entries are 0, 1, or 2
and whose row and column sums are at most 2, and that the 
sum of all the entries in $T$ is $2k$.  Then $T = \sigma + \tau$ for some
$k$-assignments $\sigma$ and $\tau$.
\end{lemma}
{\em Proof:\/}\ 
We may assume there are no 2's in $T$, since if there is a 2 we know
that both $\sigma$ and $\tau$ must have a 1 there, and only 0's everywhere else in its
row and column.  So we assume $T$ is a 0-1 matrix whose row and column sums are
at most 2, such that the sum of all entries is $2s$ for some $s\leq k$,
and we want to find two $s$-assignments $\sigma$ and $\tau$ such that $\sigma+\tau=T$.
Identify $T$ with a graph with vertices at each 1 of $T$, and edges between
any two 1's belonging to the same row or column.
Clearly every vertex of $T$ has degree $\le 2$, so every component of $T$
is a chain or a cycle.
The vertices in each component can be
alternately assigned to $\sigma$ and $\tau$.  If there is an odd component (which must be
a chain) there must
be another odd component to balance it out (so one can have an extra $\sigma$
vertex, and the other can have an extra $\tau$ vertex), since the sum of all
entries in $T$ is even. $\Box$\medskip

\begin{lemma}\label{lem:1}
Suppose that $Y$ is a $k$-reduced matrix and
$$
S=\left[\begin{array}{rr}
   a & b\\
   c & d\\
\end{array}\right]
$$
is a submatrix of $Y$.  If $a$ and $d$ each participate in
a minimum $k$-assignment of $Y$, then $a+d \le b+c$.  If also
$a+d = b+c$, then $b$ and $c$ also each participate in 
a minimum $k$-assignment.  These statements also hold for
$a$ and $d$ switched with $b$ and $c$.
\end{lemma}
{\em Proof:\/}\  
Suppose that $a+d > b+c$. Let $\sigma_1$ and $\tau_1$ be
minimum assignments passing through $a$ and $d$ respectively.
Form the matrix $T_1=\sigma_1+\tau_1$.  Then form
$T$ from $T_1$ by subtracting 1 at the positions of $a$ and $d$ 
and adding 1 at the positions of $b$ and $c$.
The hypotheses of Lemma~\ref{lem:0} apply to $T$ so that 
$T=\sigma+\tau$ for some $k$-assignments $\sigma$ and $\tau$.
But since $a+d >b+c$ we must have 
$$\sigma_1 \cdot Y + \tau_1 \cdot Y=T_1 \cdot Y > T \cdot Y=\sigma\cdot Y+\tau \cdot Y,$$
contradicting the minimality of the assignments $\sigma_1$ and $\tau_1$.
Thus $a+d \le b+c$.  

Now suppose that $b+c = a+d$.  Then the same construction yields
$$\sigma_1 \cdot Y + \tau_1 \cdot Y=T_1 \cdot Y = T \cdot Y=\sigma\cdot Y+\tau \cdot Y,$$
so that both $\sigma$ and $\tau$ are minimum $k$-assignments, and
at least one includes $b$ and at least one includes $c$.

This proof obviously also holds with $a$ and $d$ switched with $b$ and $c$.
$\Box$\medskip

\begin{proposition}\label{pro:3}
Suppose that $Y=(y_{ij})$ is a $k$-reduced $m \times n$ matrix.
Then there exist $\lambda_1,\dots,\lambda_m$ and $\mu_1,\dots,\mu_n$
such that
\begin{equation}\label{eq:3}
y_{ij}=\max(0,\lambda_i+\mu_j)
\end{equation}
and such that $y_{ij}$ participates in a minimum $k$-assignment precisely
when $y_{ij}=\mu_i+\lambda_j$.
\end{proposition}
{\em Proof:\/}\
Let $d=y_{tu}$ denote the largest entry in $Y$.
Take $\lambda_i$ to be the $i^{th}$ entry in the column of $d$, so
that $\lambda_i= y_{iu}$. Let
$\mu_j$ to be the $j^{th}$ entry in the row of $d$, decreased by $d$,
so $\mu_j=y_{tj}-d$.

First we prove (\ref{eq:3}).

When $y_{ij}$ is in the row or column of $d$, then we have
$y_{ij}=\lambda_i + \mu_j$, so (\ref{eq:3}) is immediate.

Suppose that $a$ is in neither the row of $d$ nor the column of $d$.
Let
$$
S=\left[\begin{array}{rr}
   a & b\\
   c & d\\
\end{array}\right]
$$
be the submatrix of $Y$ containing the rows and columns of $a$ and $d$
(where the order of the rows or columns in $S$ may be opposite to 
the order they occur in $Y$).
Then we just need to show that $a=\max(0,b+c-d)$.

First suppose that $b+c-d>0$.  Then $b,c,d$ must all be positive, because
$d$ is maximum for the whole matrix.
Since $b,c$ are both positive, they must both participate in a minimum
$k$-assignment, so $b+c \le a+d$, by Lemma~\ref{lem:1}.
But then, $a\ge b+c-d$, so $a$ is positive.  Then
$a$ and $d$ both participate in a minimum $k$-assignment,
which implies $a+d \le b+c$, so $a=b+c-d$.  

Next, suppose that $b+c-d \le 0$.  If $a>0$, then $a$ and $d$ both
participate in a minimum $k$-assignment, so $d<a+d \le b+c$,
a contradiction.  Thus $a=0$.

Now we show that $y_{ij}=\lambda_i+\mu_j$ exactly when
$y_{ij}$ participates in a minimum $k$-assignment.

Recall that for any $y_{ij}$ in either the row or column of $d$, we
have $y_{ij}=\lambda_i+\mu_j$.  So we need to show that all such
entries participate in a minimum $k$-assignment.  Let $b$ be any entry
in the column of $d$.  We already know that positive entries must
participate in a minimum $k$-assignment, so we assume that $b=0$.  If
$d=0$, then our whole matrix is zero and our result is trivial, so we
may assume that $d>0$ and therefore participates in a minimum
$k$-assignment.  If the minimum assignment using $d$ does not use the
row of $b$, we can replace $d$ by $b$ and obtain a smaller assignment,
a contradiction.  So we can conclude that the minimum assignment using
$d$ also uses an element $a$ from the row of $b$.  Form the $2\times
2$ submatrix $S$ containing $a$ and $d$ as above.  Since $d$ is the
largest entry, $c\le d$, and we also have $b=0 \le a$.  Thus we can
exchange $a$ and $d$ for $b$ and $c$, to obtain a minimum
$k$-assignment in which $b$ participates.  In the same way, we see that
any entry in the row of $d$ must participate in a minimum
$k$-assignment.

Finally, consider an element $a$ that is neither in the row nor the
column of $d$ and form the $2\times 2$ submatrix $S$ containing $a$ and $d$
as above.
We must show that $a$ participates in a minimum assignment
exactly when $a=b+c-d$.  This is certainly true if $a>0$.  So,
let us assume that $a=0$.  Also, since both $b$ and $c$ are
in a row or a column of $d$, both participate in a minimum
assignment, so that $b+c \le a+d$.

Now suppose that $a$ participates in a minimum $k$-assignment.  Then
$a+d \le b+c$, so $a=b+c-d$, as required.

Conversely, suppose that $a=b+c-d$. Then, since $b$ and $c$ participate in 
minimum $k$-assignments, Lemma~\ref{lem:1} shows that $a$ also
participates.  $\Box$\medskip

{\em Proof of Theorem~\ref{th:1}:\/}\ 
Now we are ready to prove Theorem~\ref{th:1} by showing that any
$k$-reduced matrix
$Y$ in $P_\sigma$ is a nonnegative linear combination of a suitable collection
of matrices $V_{IJ}$ from $P_\sigma$.
Without loss of generality we can assume that the $\lambda$'s and $\mu$'s are
weakly increasing.
In this case the rows and columns of $Y$ are also weakly increasing.  

If the matrix $Y$ is zero, there is nothing to prove.  
Otherwise, since every nonzero entry in $Y$ participates in
a minimum $k$-assignment, we know that the minimum $k$-assignment
is nonzero.  Hence 
there is
at least one nonzero entry among $y_{1,k},y_{2,k-1},\dots,y_{k,1}$.
In particular, there is a pair $(i,j)$ such that 
$y_{ij}>0$ and $i+j \le k+1$.  Now select such a pair $(i,j)$ to be minimal in 
the sense that 
if $i \ne 1$ then $y_{i-1,j}=0$ and if $j\ne 1$ then $y_{i,j-1}=0$.
Let $I=\{1,\dots,i-1\}$ and let $J=\{1,\dots,j-1\}$.  We will show
that $V_{IJ}$ is in $P_\sigma$ and that
$Y-y_{ij}V_{IJ}$ is again in $P_\sigma$ and still $k$-reduced.

Suppose that $1\le i'<i$ and $1\le j'<j$.
Then, since 
$$  \max(0,\lambda_{i'}+\mu_j)=y_{i'j}=0<y_{ij}=\max(0,\lambda_i+\mu_j),  $$
we know that $\lambda_{i'}<\lambda_i$.  Similarly
$\mu_{j'} < \mu_j$.
Since $y_{i'j}=\max(0,\lambda_{i'}+\mu_j)=0$, we have
$\lambda_{i'}+\mu_j \le 0$.  But then
$$ \lambda_{i'}+\mu_{j'} <  \lambda_{i'}+\mu_{j}\le 0 .$$
If follows from Proposition~\ref{pro:3} that
none of the matrix entries $y_{i'j'}$ 
with $i' < i$ and
$j' < j$ can participate in a minimum $k$-assignment. 

To see that $V_{IJ}\in P_{\sigma}$, first note that any minimum
$k$-assignment of $Y$ must use all of the first $i-1$ rows.  If not,
since $i \le k$, there is some $i_1>i$ such that row $i_1$
participates.  Then, since not all of the first $i-1$ rows are used,
we can replace the entry of the assignment in row $i_1$ with the entry
in the same column of row $i_0$, for some $i_0<i$, to get an
assignment with a value no larger.  The entry being replaced could not
come from a column preceding $j$, by the discussion above, so it must
be positive.  But also by the discussion above, $\lambda_{i_0} <
\lambda_{i_1}$.  Then the new assignment would be strictly smaller, a
contradiction.  Thus, any minimum assignment uses all of the first
$i-1$ rows and all of the first $j-1$ columns, and does not use any
entry which is in both the first $i-1$ rows and the first $j-1$
columns.  So, if $\tau$ is any matrix representing a minimum
$k$-assignment of $Y$, we must have $V_{IJ}\cdot \tau=k-|I|-|J|$.  In
particular, $V_{IJ}\cdot \sigma=k-|I|-|J|$, so $V_{IJ}$ is in
$P_\sigma$ as claimed.

The preceding argument shows that if we replace $Y$ by $Y'=Y-tV_{IJ}$, for
any $t$ satisfying $0\le t\le y_{ij}$, the effect on any minimum $k$-assignment is to subtract
$(k-|I|-|J|)t$ from its value.
Thus, all assignments $\tau$ that are minimum
for $Y$ will agree on $Y'$.  We now show that each of these assignments
$\tau$ is minimum for $Y'$ as well.  Assume not.
Then we could find $t_1$ and $t_2$ such that  $ 0\le t_1 < t_2\le y_{ij}$ and
$k$-assignments $\tau$ and $\phi$ such that
$\tau$ is minimum for $(Y-tV_{IJ})$ when $ t \leq t_1$ but not
when $t_1 < t \leq t_2$, and
$\phi$ is minimum for $Y-tV_{IJ}$ when $t_1 \le t \le t_2$.
Thus, we would have both $\phi$ and $\tau$ minimum assignments
for $Y-t_1V_{IJ}$.  But then our preceding argument applied to
$Y-t_1V_{IJ}$ tells us that
$\phi$ and $\tau$ must agree on $Y-tV_{IJ}$ for $t_1 \le t \le t_2$,
a contradiction.

Now we let $Y'=Y-y_{ij}V_{IJ}$ and observe that $Y'$ is $k$-reduced.
Indeed, all the assignments $\tau$ that were minimum for
$Y$ are also minimum for $Y'$.
Thus, any element that participated in a minimum $k$-assignment for 
$Y$ will also participate in a minimum $k$-assignment for $Y'$.
Also, the replacement of $Y$ by $Y'$ creates no new nonzero elements,
so the matrix $Y'$ will be $k$-reduced. 

Since $\sigma$ in particular is a minimum assignment for $Y$,
$\sigma$ will also be minimum for $Y'$, so that 
$Y'$ is also in $P_\sigma$.

Note that $Y'$ has at least one more zero entry than $Y$,
namely the entry at $(i,j)$. 

We can continue removing multiples of submatrices $V_{IJ}$, each
time producing a matrix with at least one more zero entry.  Thus,
we eventually reach a matrix that is all zero.  In effect,
we have expressed $Y$ as a nonnegative linear combination of the
generators as required.  $\Box$\medskip

We have shown that every element of $P_\sigma$ is a nonnegative linear 
combination of certain $e_{ij}$'s and $V_{IJ}$'s in $P_\sigma$.  We remark
that the $e_{ij}$'s generate extreme rays of $P_\sigma$ but the 
$V_{IJ}$'s in general do not.  It is not hard to show that the
$V_{IJ}$'s which do generate extreme rays of $P_\sigma$ are
those with $|I|+|J|=k-1$.

We now observe that Theorem~\ref{th:1} allows us to make some conclusions
about the rational function $E_k(A)$.  A simplicial cone in
a decomposition of $P_\sigma$ has some generators of the form $e_{ij}$ and
some of the form $V_{IJ}$.  For the generators of the form $e_{ij}$, we know
that $\sigma\cdot e_{ij}=0$. Also the dot products $e_{ij} \cdot V$
in the denominator cancel with factors in
the initial product of $a_{ij}$'s in (\ref{eq:2}).
Thus, the denominator of the integral over a simplicial cone is
a product of terms of the form $A \cdot V_{IJ}$.  Finally, we can
conclude:

\begin{theorem}\label{th:genform}
The expected minimum $k$-assignment of a random exponential matrix with
rate matrix $A$ is a rational function of the entries of $A$. The
denominator of the rational function is a product of sums, each
being the sum of all entries in a submatrix of $A$ omitting $i$ rows
and $j$ columns, where $i+j<k$.
\end{theorem}

\begin{example}\label{ex:arb2}
The minimum 2-assignment of a random $2\times2$ exponential matrix with
rate matrix $A=(a_{ij})$ is
\begin{eqnarray*}
\frac1{a_{11}+a_{12}}&+&\frac1{a_{21}+a_{22}}
+\frac{a_{11}a_{21}}{(a_{11}+a_{12})(a_{21}+a_{22})(a_{12}+a_{22})}\\
&+& \frac{a_{12}a_{22}}{(a_{11}+a_{12})(a_{21}+a_{22})(a_{11}+a_{21})}
\end{eqnarray*}
\end{example}
The above formula is easily computed by the method we will sketch
in Section~\ref{sec:comp}.

In the case of a rank 1 rate matrix $A=(r_ic_j)$,
Theorem~\ref{th:genform} says that $E_k(A)$ is a rational function
of the $r_i$'s and $c_j$'s and its denominator is
a product of sums of subsets of the $r_i$'s omitting fewer than $k$
of the $r_i$'s and sums of subsets of the $c_j$'s omitting fewer than $k$
of the $c_j$'s.
Since this is a consequence of Conjecture~\ref{conj:r1}, it lends some
support to the conjecture.

\subsection{The nesting lemma}\label{sec:nest}

A real nonnegative $m \times n$ matrix $X$ has minimum $k$-assignments for
each $k\le m$.  Generically there
is only one of each but in some cases there are many minimum assignments
of various sizes.  It helps to know how these are related.

The following lemma is fundamental.  Other proofs probably exist
but we include ours here for completeness.

\begin{lemma}\label{lem:nest}
Let $k_1$ and $k_2$ be two integers, with $k_1 \le k_2 \le m$.

Suppose that $M_1$ is a $k_1 \times k_1$ submatrix of $X$ that
contains a minimum $k_1$-assignment of $X$.  Then there exists a 
$k_2 \times k_2$ submatrix $M_2$ containing $M_1$ such that $M_2$
contains a minimum $k_2$-assignment of $X$.

Suppose that $M_2$ is a $k_2 \times k_2$ submatrix of $X$ that
contains a minimum $k_2$-assignment of $X$.  Then there exists a 
$k_1 \times k_1$ submatrix $M_1$ contained in $M_2$ such that $M_1$
contains a minimum $k_1$-assignment of $X$.
\end{lemma}
{\em Proof:\/}\ 
If $k_1=k_2$, there is nothing to prove. 
So assume $k_1 < k_2$ and 
fix a minimum $k_1$-assignment and a minimum $k_2$-assignment.

Let $G$ be the graph on $k_1$ red vertices (representing the entries
of the $k_1$-assignment) and $k_2$ blue vertices (representing the
entries of the $k_2$-assignment), with edges between two vertices if
the corresponding entries belong to the same row or column.  (If the
assignments share an entry then we have a red vertex and a blue vertex
with two edges between them comprising a component which is a cycle of
length 2.)

Then $G$ is bipartite and no vertex of $G$ has degree more than 2.
Thus, every component of $G$ is a cycle or a chain in which the red
and blue vertices alternate.

Suppose some component of $G$ has $m_1$ red vertices and $m_2$ blue vertices,
and all of its red vertices have degree 2.
Such a component is either a cycle, so that $m_1=m_2$,
or a chain with blue vertices at each end, so that $m_1+1=m_2$.
In either case we have $m_1+1 \ge m_2$.

Consider the submatrix $M_1$ spanned by the associated $m_1$ entries of the
$k_1$-assignment and the submatrix $M_2$ spanned by the associated $m_2$ entries
of the $k_2$-assignment.  The component condition translates
into the condition that the $M_1$ is contained in $M_2$.
Also the remaining entries of the two assignments
comprise a minimum $(k_1-m_1)$-assignment and a minimum 
$(k_2-m_2)$-assignment of the submatrix of $X$ complementary to $M_2$.
Since $k_1 < k_2$ and $m_1+1 \ge m_2$, we have $k_1-m_1 \le k_2-m_2$, 
so the lemma follows by induction applied to the
complementary submatrix.

Now assume that every component of $G$ has a red vertex of degree 1 or
less.  Such components are chains with one endpoint red.  When the
other endpoint is red, there are more red than blue vertices in the
component.  When the other endpoint is blue, the number of vertices of
both colors is equal.  In particular the number of red vertices is
always at least as great as the number of blue vertices.  It follows
that, all together, there are at most $k_1$ blue vertices that are
connected to some red vertex.  Thus we can select $k_2-k_1$ entries
of the $k_2$-assignment that do not share any row or column with the
$k_1$-assignment.  We can now consider three sets of matrix entries,
the set $S_1$ of $k_1$ entries of the $k_1$-assignment, the set $S_2$,
just selected, of $k_2-k_1$ of entries of the $k_2$-assignment and the
set $S_3$ of the remaining $k_1$ entries of the $k_2$-assignment.
Then the sets $S_3$ and $S_1$ are both in the submatrix complementary
to that determined by $S_2$.  Moreover the sums of the entries in
$S_1$ and $S_3$ must be equal.  For, if the sum of the $S_1$ entries
were greater than that of $S_3$, then $S_3$ would be a smaller
$k_1$-assignment than $S_1$, a contradiction.  But if the sum of the
$S_1$ entries were smaller than the sum of $S_3$ entries, then $S_1
\cup S_2$ would be give a smaller $k_2$-assignment than $S_2 \cup
S_3$, which was a minimum $k_2$-assignment.

Thus $S_3$ is a minimum $k_1$-assignment contained in our minimum
$k_2$-assignment and $S_1 \cup S_2$ is a minimum $k_2$-assignment containing our
minimum $k_1$-assignment, which proves our lemma. $\Box$

\subsection{A computational consequence of the nesting lemma}\label{sec:aux}

Lemma~\ref{lem:nest} suggests a way to compute the minimum
$k$-assignment of a matrix $X$ when $k<n$.
We proceed by finding the
minimum $k$-assignments for $X$ for $k=1,2,\dots,$ one at a time.
Suppose we have found a minimum $(k-1)$-assignment whose rows and
columns determine a $(k-1) \times (k-1)$ submatrix $M$ of $X$.  Then, when we
search for a minimum $k$-assignment, we know that we can restrict our
search to minimum $k$-assignments in one of the submatrices
of $X$ obtained by appending a single new row and new column to $M$.

There are some simple properties that such an extension must have.
Suppose the submatrix of a minimum $k$-assignment uses $M$ together
with a new row $i$ and a new column $j$.  Also, suppose that the minimum
$k$-assignment uses an entry from column $j$ that is in row $i'$ of
$M$.  Then this entry in column $j$ must be the smallest entry in the
part of row $i'$ outside of $M$.  Similarly, if the minimum
$k$-assignment uses an entry from row $i$ that is in column $j'$ of
$M$, then this entry in row $i$ must be the smallest entry in the part
of column $j'$ outside $M$.  Finally, if the minimum $k$-assignment
uses the entry $x_{ij}$, then this entry must be minimum in the
submatrix of $X$ complementary to $M$.

The preceding discussion shows that the following strategy will work
to construct the minimum $k$-assignment once we have found the
minimum $(k-1)$-assignment.
We define the $k \times k$ \emph{auxiliary matrix} $\mbox{Aux}_M(X)$ by
appending to $M$ a new row and column as follows.  To each row $i$ of
$M$ we append a new entry which is the minimum of the entries
in row $i$ of $X$ that are outside of $M$.  To each column $j$ of $M$
we append a new entry which is the minimum of the entries in column
$j$ of $X$ that are outside of $M$.  At the intersection of the new
row and new column we place the minimum entry of the submatrix of $X$
complementary to $M$.

We now find the minimum $k$-assignment of $\mbox{Aux}_M(X)$.  This
assignment will use an entry in the last row and an entry in the
last column of $\mbox{Aux}_M(X)$, which can be the same entry
if the assignment uses the entry in the last row and column.
Each of these entries is a copy of some entry of $X$,
which then tells us which row and column of $X$ need to be appended
to $M$ to obtain the submatrix of the minimum $k$-assignment of $X$.

\subsection{A consequence for the expected contribution}

The discussion of the auxiliary matrix in the preceding
section can be formalized to prove an interesting property
of the expected contribution.

We will use the following simple facts about independent
exponential random variables.

\begin{proposition}\label{pro:prob}
Suppose that $a_1,\dots,a_m$ are positive real numbers,
and that $x_1,\dots,x_m$ are independent random variables 
with $x_i$ chosen
from the exponential distribution with rate $a_i$.

Let $x$ denote the random variable $\min_i x_i$.  Then $x$ is distributed as an
exponential random variable of rate $a_1+\cdots+a_m$.

Let $W$ be the discrete random variable whose value is the least $i$
for which $x=x_i$. Then $W$ and $x$ are independent random
variables, and the probability that $W=i$ is $a_i/(a_1+\cdots+a_m)$.
\end{proposition}
{\em Proof:\/}\  
Let $c$ be a positive real number.
Then the probability that $x\ge c$ is
\begin{eqnarray*}
& & a_1 \cdots a_m \int_{x_1=c}^\infty \cdots \int_{x_m=c}^\infty  e^{-\sum_{j=1}^m a_jx_j} dx_1\cdots dx_m. \\
&=& a_1 \cdots a_m e^{-c\sum_{j=1}^m a_j} 
 \int_{u_1=0}^\infty \cdots \int_{u_m=0}^\infty e^{-\sum_{j=1}^m a_ju_j} du_1\cdots du_m. \\
&=& e^{-c\sum_{j=1}^m a_j}, 
\end{eqnarray*}
where we have made the substitution $x_i=u_i+c$ on the second line.
Taking the derivative with respect to $c$ we see that $x$ is distributed
as an exponential random variable of rate $a_1+\cdots+a_m$.

A similar substitution yields
the part of the integral corresponding the event that $x=x_1$, as follows:
\begin{eqnarray*}
& & a_1 \cdots a_m \int_{x_1=c}^\infty \int_{x_2=x_1}^\infty \cdots \int_{x_m=x_1}^\infty  e^{-\sum_{j=1}^m a_jx_j} dx_1\cdots dx_m. \\
&=& a_1 \cdots a_m  \int_{u_1=0}^\infty \cdots \int_{u_m=0}^\infty 
e^{-(a_1(c+u_1)+a_2(c+u_1+u_2) +\cdots a_m(c+u_1+u_m))} du_1\cdots du_m  \\
&=& \frac{a_1}{\sum_{j=1}^m a_j} e^{-c\sum_{j=1}^m a_j}.
\end{eqnarray*}
Thus, the conditional probability that $x=x_i$, given $x\geq c$, is
$a_i/(a_1+\cdots+a_m)$. Since this is true for all $c$, the event
$x=x_i$ is independent of the random variable $x$. $\Box$\medskip

Now let $X$ be any nonnegative $m\times n$ matrix and
associate to $X$ the $k\times k$ matrix $Y=\mbox{Aux}_M(X)$
where $M$ is the upper left $(k-1)\times(k-1)$ submatrix of $X$.
For any $t<k$, let $M_t$ denote the upper left $t\times t$ submatrix of $X$.
By abuse of notation we will also let $M_t$ denote the upper left
$t\times t$ submatrix of $Y$, since they are identical.

\begin{lemma}\label{lem:auxt}
For any $t< k$, let $\sigma$ be a $t$-assignment of $M_t$.
Then $\sigma$ is a minimum $t$-assignment for $X$ if and only if it is
a minimum $t$-assignment for $Y$.
\end{lemma}
{\em Proof:\/}\ 
Let $\sigma$ be a $t$-assignment of $M_t$.

If $t=1$ it is clear that the lemma holds. We proceed by induction on $t$.

Suppose that $\sigma$ is a minimum $t$-assignment of $X$.
By induction on $t$, there is a minimum $(t-1)$-assignment
of $Y$ which lies in $M_t$, so by Lemma~\ref{lem:nest}, there is
some minimum $t$-assignment $\tau$ of $Y$ which uses at most one
column and one row outside of $M_t$.  
By definition, each entry of $Y$ outside of $M_t$ equals some entry
of $X$ outside of $M_t$.
If we replace each of the entries of $\tau$ outside of $M_t$ by
their equivalent entries in $X$, then we get a $t$-assignment $\tau'$
of $X$.  Furthermore, $\tau'\cdot X \ge \sigma \cdot X$, since
$\sigma$ is minimum.  But $\tau'\cdot X = \tau\cdot Y$ and $\sigma\cdot X
= \sigma\cdot Y$, so $\tau \cdot Y \ge \sigma \cdot Y$.
Therefore, $\sigma$ is a minimum $t$-assignment of $Y$.

Now suppose that $\sigma$ is a minimum $t$-assignment of $Y$.
By induction on $t$, there is a minimum $(t-1)$-assignment
of $X$ which lies in $M_t$, so by Lemma~\ref{lem:nest}, there is
some minimum $t$-assignment $\tau$ of $X$ which uses at most
one column and one row outside of $M_t$.
We can replace each entry $x_{ij}$ of $\tau$
by the entry $y_{\min(i,k),\min(j,k)}$ to get a $t$-assignment
$\tau'$ of $Y$.  Furthermore, $\tau'\cdot Y\ge \sigma \cdot Y$,
since $\sigma$ is minimum.  But $\tau\cdot X \ge \tau'\cdot Y$
and $\sigma\cdot X = \sigma\cdot Y$, so
$\tau \cdot X \ge \sigma \cdot X$.   Therefore,
$\sigma$ is a minimum $t$-assignment of $X$.
$\Box$\medskip

\begin{lemma}\label{lem:aux}
The minimum $(k-1)$-assignment of $X$ uses its first $k-1$
rows and columns exactly when the minimum $(k-1)$-assignment of $Y$
uses its first $k-1$ rows and columns.  In this case,
the value of the minimum $k$-assignment of $X$ is the 
same as the value of the minimum $k$-assignment of $Y$.

Moreover 
in this case there is a minimum $k$-assignment $\tau$ of $X$
and a minimum $k$-assignment $\tau'$ of $Y$ which correspond entry by entry;
{\em i.e.,\/}\ the entries of $\tau$ in $M$ are located in the same
positions as the entries of $\tau'$ in $M$, and for any entry $x_{ij}$ of
$\tau$ outside of $M$ there is a corresponding entry
$y_{\min(i,k),\min(j,k)}$ of $\tau'$ outside of $M$.
\end{lemma}
{\em Proof:\/}\ 
The first paragraph is the case $t=k-1$ in Lemma~\ref{lem:auxt}.
The second paragraph follows from our discussion
of auxiliary matrix in the preceding section. 
$\Box$\medskip

We will use the more detailed statement in the second paragraph 
to prove Lemma~\ref{lem:flag1}
in Section~\ref{sec:flag}.  

Now let $X$ be a random exponential $m\times n$ matrix
with rate matrix $A$, and associate to $A$ the $k \times k$ matrix
$B=(b_{ij})$ with 
$$  b_{ij}=a_{ij}\quad \hbox{when $1 \le i,j \le k-1$}  $$
and 
$$ b_{ik}=\sum_{j' \ge k}a_{ij'}\,,\quad i=1,\dots,k-1 $$
and 
$$ b_{kj}=\sum_{i' \ge k}a_{i'j}\,,\quad j=1,\dots,k-1 $$
and
$$  b_{kk}=\sum_{i',j' \ge k} a_{i'j'}.  $$
Note that $Y=\mbox{Aux}_M(X)$ is a random exponential matrix with
rate matrix $B$, by Proposition~\ref{pro:prob}.

The statement in the first paragraph of the preceding Lemma shows that
the expected contribution of the submatrix of $X$ consisting of its first $k-1$ rows and
first $k-1$ columns to the minimum $k$-assignment of $X$ is the same
as the expected contribution 
of the submatrix of $Y$ consisting of its first $k-1$ rows and first $k-1$ columns
to the minimum $k$-assignment of $Y$.  


In the rank 1 rate matrix case, if $C(k,r_1,\dots,r_m,c_1,\dots,c_n)$
denotes the expected contribution of the submatrix of the first 
$k-1$ rows and columns when the rate matrix is $A=(r_ic_j)$, then
$$C(k,r_1,\dots,r_{k-1},r_k+\cdots+r_m,c_1,\dots,c_{k-1},c_k+\cdots+c_n)$$
is the expected contribution when the rate matrix is $B$.
These two functions must agree, so the entries of $A$ outside of the first $k-1$ rows and columns
enter into the contribution function only via the sums defining
the entries in the $k^{th}$ row and column of $B$.

Now note that Conjecture~\ref{conj:contrib} has the feature that
it predicts the same contribution in the two cases above.  So this
is a bit of evidence in favor of Conjecture~\ref{conj:contrib}.

Thus, we can summarize our discussion in the rank 1 rate matrix case
as follows:
\begin{theorem}
Let $k$, $m$ and $n$ be integers $k \le m \le n$.  Then
Conjecture~\ref{conj:contrib} holds for $k$-assignments in an $m
\times n$ matrix if and only if it holds for $k$-assignments in a $k
\times k$ matrix.
\end{theorem}

Thus, if one could prove only the cases $k=m=n$ of 
Conjecture~\ref{conj:contrib}, that would prove the
general case of that conjecture as well as Conjecture~\ref{conj:r1} and
Conjecture~\ref{conj:cs}.

\section{Rank 1 rate matrices}\label{sec:rank1}

In this section we restrict our attention to 
random exponential matrices for which the rate matrix has rank 1.
We describe how we arrived at Conjecture~\ref{conj:r1} and then give 
equivalent formulations, which provide different kinds of confirmation.
Finally, we prove a probability result about the locations of the minimum
$\ell$-assignments for $1\le \ell \le k$.

We discovered Conjecture~\ref{conj:r1} while experimenting with the
computations described in Section~\ref{sec:comp}.  Finding that
Mathematica had trouble carrying out the computation when the rate
matrix consisted of $mn$ indeterminates, we decided to try a simpler
case, with rate matrix of the form $a_{ij}=a_i$ for all $i$ and
$j$.  We noticed in this case that the answer has a surprisingly simple
form---in particular, it can be written as a linear combination of the
reciprocals of sums of the $a_i$'s.  Next we found that, when the rate
matrix has rank 1, so that $ a_{ij}=r_ic_j $, the expected value seems
to be a linear combination of terms of the form
\begin{equation}\label{eq:recipsum}
 \frac{1}{(\sum_{i \notin I} r_i)(\sum_{j \notin J} c_j)} 
\end{equation}
with $I$ a proper subset of $\{1,\dots,m\}$, and $J$ a proper subset
of $\{1,\dots,n\}$.
It is easily shown that the rational functions
$\frac{1}{(\sum_{i \notin I} r_i)(\sum_{j \notin J} c_j)} $
are linearly independent over the real numbers.  Thus
the coefficients in such a linear combination are uniquely determined.

Making the assumption that the expected value is indeed a linear
combination of terms of the form~(\ref{eq:recipsum}), we arrived
at Conjecture~\ref{conj:r1} by considering certain limiting
conditions on the expected value.  We will describe these limiting 
conditions in Section~\ref{sec:extra}.

\subsection{Equivalent formulations of Conjecture~\ref{conj:r1}}
\label{sec:equiv}

From now on we use the shorthand notation $[m]$ for the set $\{1,\dots,m\}$.

Let us introduce the notation 
\begin{equation}\label{eq:formula}
 F(k,r,c)=
\sum_{I,J}
  (-1)^{ k - 1 - |I| - |J| } \cdot 
   \binom{m + n - 1 - |I| - |J|}{k - 1 - |I| - |J|}
      \frac{1}{( \sum_{i \notin I}r_i) \cdot ( \sum_{j \notin J} c_j )}
\end{equation}
for the formula given in Conjecture~\ref{conj:r1}.  Here recall
that $r=(r_1,\dots,r_m)$ is an $m$-tuple of positive real numbers,
$c=(c_1,\dots,c_n)$ is an $n$-tuple of positive real numbers, and 
the sum is over proper subsets $I\subsetneq [m]$ and $J\subsetneq[n]$.  
The binomial coefficient enforces the condition
$|I|+|J|<k$.
In what follows we will often not mention
such constraints explicitly.

In this section we derive alternative ways to write
(\ref{eq:formula}) and conclude that $F(k,r,c)$ is positive,
Conjecture~\ref{conj:r1} implies Conjecture~\ref{conj:cs}, and
Conjecture~\ref{conj:contrib} implies Conjecture~\ref{conj:r1}.

Note that $F(k,r,c)$ can be written more succinctly as
\begin{equation}\label{eq:formula1}
 F(k,r,c)=
\sum_{I,J}
   \binom{k-1-m-n}{k - 1 - |I| - |J|}
      \frac{1}{( \sum_{i \notin I}r_i)  ( \sum_{j \notin J} c_j )},
\end{equation}
using binomial coefficients with negative numerator.

\begin{proposition}\label{lem:incl-excl}
\begin{equation}\label{eq:formula2}
F(k,r,c)=\sum_{ |I'|+|J'| < k, I \subseteq I', J \subseteq J' }
          (-1)^{|I'|-|I|+|J'|-|J|} \frac{1}{ (\sum_{i \notin I} r_i ) \cdot
              ( \sum_{j \notin J} c_j ) }.
\end{equation}
Here $I'$ and $I$ are proper subsets of $[m]$ and $J'$ and $J$ are
proper subsets of $[n]$.
\end{proposition}
{\em Proof:\/}\ 
Comparison with (\ref{eq:formula1}) shows that, for
a fixed $I$ and $J$ with $|I|+|J|<k$, we need to evaluate the sum
$$\sum_{I \subseteq I',J \subseteq J',|I'|+|J'|<k } (-1)^{|I'|-|I|+|J'|-|J|}.  $$
If we denote by $i$ and $j$ the cardinalities of $I$ and $J$ and by $t$ and
$u$ the cardinalities of $I'$ and $J'$, we can rewrite this sum as 

\begin{eqnarray*}
& &\sum_{t \ge i,u\ge j, t+u<k} (-1)^{t-i+u-j}\binom{m-i}{t-i}\binom{n-j}{u-j}  \\
&=&\sum_{t \ge 0,u\ge 0, t+u<k-i-j} (-1)^{t+u}\binom{m-i}{t}\binom{n-j}{u}  \\
&=&\sum_{l=0}^{k-1-i-j}(-1)^l \sum_{t \ge 0,u\ge 0, t+u=l}
\binom{m-i}{t}\binom{n-j}{u}  \\
&=& \sum_{l=0}^{k-1-i-j} (-1)^l\binom{m+n-i-j}{l} \\
&=& \sum_{l=0}^{k-1-i-j} \binom{i+j-m-n+l-1}{l} \\
&=& \sum_{l=0}^{k-1-i-j} \left( \binom{i+j-m-n+l}{l}- \binom{i+j-m-n+l-1}{l-1}\right)\\
&=& \binom{k-1-m-n}{k-1-i-j}  \\
\end{eqnarray*}
which agrees with (\ref{eq:formula1}). $\Box$\medskip

We can rewrite (\ref{eq:formula2}) as a double sum with the inner
sum over $I$ and $J$ and the outer sum over $I'$ and $J'$.  Then, for fixed
$I'$ and $J'$, the inner sum factors as
\begin{equation}\label{eq:inex}
 \left( \sum_{I \subseteq I'} (-1)^{|I'|-|I|}\frac{1}{\sum_{i \notin I}r_i}\right)
 \left( \sum_{J \subseteq J'} (-1)^{|J'|-|J|}\frac{1}{\sum_{j \notin J}c_j}\right)
\end{equation}
Now we show that each factor has an interesting probabilistic 
interpretation.  

Suppose that an urn contains $m$ balls labeled $1,2,\dots,m$ and
for each $i$, ball $i$ has weight $r_i$.
We select balls one at a time without replacement,
at each time selecting a ball with probability proportional to the weights
of those balls still in the urn.  
Let $\Pr(r,I')$ denote the probability that the set of balls in $I'$ are the
first $t$ balls to be chosen, where $t$ is the cardinality of $I'$.  Then
\begin{equation}\label{eq:prRI}
\Pr(r,I')=\sum_\pi \prod_{i=1}^t \frac{ r_{\pi_i}}{R-\sum_{j=1}^{i-1} r_{\pi_i}}  
\end{equation} 
where $R=\sum_{i=1}^m r_i$ and
the outer sum is over all $t!$ orderings $(\pi_1,\dots,\pi_t)$ of $I'$.

We can calculate $\Pr(r,I')$ in a different way as follows.  Suppose
we draw all $m$ balls from the urn.  If we fix any subset $U$ of balls, then
the probability that a particular ball $u$ from $U$ is chosen before
any other ball from $U$ is the weight of $u$
divided by the sum of the weights of the balls in $U$.

Now, for $i\in I'$,
let $E_i$ denote the event that the first time a ball is drawn from the
set consisting of $i$ together with the complement of $I'$, 
the ball chosen is from the complement of $I'$.
Then $E_i$ has probability
$$  \frac{\sum_{j\notin I'} r_j}{r_i + \sum_{j \notin I'} r_j}.  $$
In order for our set $I'$ to be the set of the first $t$ balls chosen, it is
necessary and sufficient that none of the events $E_i$ occur. 
For any subset $I$ of $I'$ the
probability that all of the events $E_i$, $i \in I$, occur is
$$  \frac{\sum_{j \notin I'} r_j}{\sum_{i \notin (I'- I)} r_i}  $$
So, by the Inclusion-Exclusion Principle, 
\begin{equation}\label{eq:incex}
\Pr(r,I') = \sum_{I \subseteq I'} (-1)^{|I'|-|I|}\frac{\sum_{i \notin I'}r_i}{\sum_{i \notin I} r_i} 
\end{equation}
which is  $(\sum_{i \notin I'}r_i)$ times the first factor in 
(\ref{eq:inex}).  The analogous result holds for the second factor
in (\ref{eq:inex}).   We conclude that
\begin{proposition}\label{lem:probs}
\begin{equation}\label{eq:pr1}
F(k,r,c)=
\sum_{I,J,|I|+|J| < k}
   \frac{\Pr(r,I)\Pr(c,J)}{(\sum_{i \notin I}r_i)(\sum_{j \notin J}c_j)},
\end{equation}
where the sum is over proper subsets $I$ of $[m]$ and 
$J$ of $[n]$. 
Hence $F(k,r,c)$ is always positive.
\end{proposition}
$\Box$\medskip

Now we rewrite $F(k,r,c)$ to show that Conjecture~\ref{conj:r1}
implies Conjecture~\ref{conj:cs} and Conjecture~\ref{conj:contrib}
implies Conjecture~\ref{conj:r1}.

Let $\Pr(r,(i_1,\dots,i_t))$ denote the probability
that the first $t$ selections from our urn are $i_1,\dots,i_t$ in that order.
Then, from (\ref{eq:pr1}) and (\ref{eq:prRI}), we get
\begin{equation}\label{eq:pr2}
F(k,r,c)=
\sum_{t,u\ge 0, t+u < k}\sum_{i,j}
   \frac{\Pr(r,(i_1,\dots,i_t))\Pr(c,(j_1,\dots,j_u))}
{(R-\sum_{s=1}^tr_{i_s} )
   (C-\sum_{s=1}^u c_{j_s})},
\end{equation}
where the sum is over all sequences $i=(i_1,\dots,i_t)$ of distinct integers
in $[m]$ and $j=(j_1,\dots,j_u)$ of distinct integers
in $[n]$.  But, certainly 
$$\Pr(r,(i_1,\dots,i_t))=\sum_i \Pr(r,(i_1,\dots,i_t,\dots,i_{k-1})) $$
where the sum is over all extensions of $(i_1,\dots,i_t)$ to a 
$(k-1)$-long sequence $i$ of distinct integers in $[m]$.
Thus, we can rewrite (\ref{eq:pr2}) and obtain:
\begin{proposition}
\begin{equation}\label{eq:pr3}
F(k,r,c) = 
\sum_{i,j}\sum_{t,u\ge 0, t+u < k}
   \frac{\Pr(r,(i_1,\dots,i_{k-1}))\Pr(c,(j_1,\dots,j_{k-1}))}
{(R-\sum_{s=1}^tr_{i_s} )
   (C-\sum_{s=1}^u c_{j_s})},
\end{equation}
where the outer sum in (\ref{eq:pr3}) is over pairs of ordered
sequences of $k-1$ distinct integers from $[m]$ and $[n]$.
\end{proposition}

Note that each term in the above
sum corresponds to a flag of submatrices of sizes $1\times 1,\dots,
(k-1)\times (k-1)$.
In this form, specializing to the case that all the
$r$'s and $c$'s are 1, it is easy to see that Conjecture~\ref{conj:r1} implies
Conjecture~\ref{conj:cs}.

We now group the terms in the outer sum according to the
(unordered) sets $I=\{i_1,\dots,i_{k-1} \}$ and $J=\{j_1,\dots,j_{k-1}\}$.
It then becomes 
\begin{equation}\label{eq:pr5}
F(k,r,c)=
\sum_{I,J}\sum_{i,j}\sum_{t,u\ge 0, t+u < k}
   \frac{\Pr(r,(i_1,\dots,i_{k-1}))\Pr(c,(j_1,\dots,j_{k-1}))}
{(R-\sum_{s=1}^tr_{i_s} )
   (C-\sum_{s=1}^u c_{j_s})}
\end{equation}
where the outer sum is over sets $I$ and $J$ of size $k-1$ and the
inner sum is over permutations $(i_1,\dots,i_{k-1})$ of $I$ and
permutations $(j_1,\dots,j_{k-1})$ of $J$.  In this form we can see
that the term of the outer sum corresponding to the sets $I$ and $J$
is the expected value of the contribution of the submatrix with row
indices $I$ and column indices $J$ predicted by
Conjecture~\ref{conj:contrib}.  Since the sum of the expected
contributions of all submatrices is the expected minimum
$k$-assignment, we now see that Conjecture~\ref{conj:contrib} implies
Conjecture~\ref{conj:r1}.

Finally, for any $T \subseteq I \subseteq [m]$, let
$\Pr(r,T,I)$ denote the probability that the first $|T|$ balls drawn
from the urn comprise the set $T$ and that the first $|I|$ balls drawn
comprise the set $I$.
Then we can rewrite our formula for the expected contribution
of a submatrix with rows $I$ and columns $J$ as
\begin{equation}
\sum_{T \subseteq I,U \subseteq J , |T|+|U| < k}
\frac{   \Pr(r,T,I)\Pr(c,U,J) }
{(\sum_{t \notin T}r_t )(\sum_{u\notin U}c_u)}
\end{equation}

\subsection{Flag probabilities}\label{sec:flag}

In this section, we prove ap probability result in the special case that
the rate matrix has rank 1.
This result may be the reason that simple formulas exist
for $E_k(A)$ when the matrix $A$ has rank 1.

It is possible for a matrix to have many
minimum $k$-assignments for some $k$.  However, 
with probability 1, a random matrix $X$ has a unique
minimum $k$-assignment for each $k$.  So,
if we let $M_k$ denote a $k\times k$ submatrix of $X$ containing
a minimum $k$-assignment of $X$, then, with probability 1, $M_k$
is unique.  By Lemma~\ref{lem:nest}, the submatrices $M_k$ are
nested: $M_1 \subset M_2 \subset \cdots \subset M_k. $
We will call this the \emph{flag of submatrices of $X$}.
This flag can also be described by the list $i_1,i_2,\dots,i_k$ of
appended rows and the list $j_1,j_2,\dots,j_k$ of appended columns.
Thus, $M_l$ is the submatrix with rows
$i_1,\dots,i_l$ and columns $j_1,\dots,j_l$.

It is natural to ask for the probability that 
a random matrix has a given flag of submatrices.  We know
of no formula for this probability for general rate matrices.  However,
for rate matrices of rank 1, we can prove a simple formula
for the probability of each flag.  Moreover, this formula
will involve the probabilities $\Pr(r,(i_1,\dots,i_k))$
calculated in the previous section.

We first need the following

\begin{lemma}\label{lem:flag1}
Let $X$ be an exponential random matrix with rank 1 rate
matrix $A=(r_ic_j)$.  If
the minimum $(k-1)$-assignment of $X$ uses the first
$k-1$ rows and first $k-1$ columns, then the minimum $k$-assignment
uses an additional row and column.  The probability of row $i'$ being the
additional row is $\frac{r_{i'}}{\sum_{i=k}^m r_i}$,
the probability of column $j'$ being the additional column is
$\frac{c_{j'}}{\sum_{j=k}^n c_j}$, and these events are independent.
\end{lemma}

{\em Proof:\/}\ Let $M$ denote the upper left $(k-1)\times (k-1)$
submatrix of $X$, and $Y$ the $k\times k$ matrix
$\mbox{Aux}_M(X)$ defined in Section \ref{sec:aux}.
Then the upper left submatrix of $Y$ is identical to $M$, and, by
abuse of notation, we will let $M$ denote that submatrix of $Y$ also.
If the minimum $(k-1)$-assignment of $X$ lies in $M$,
then by Lemma~\ref{lem:auxt}, the minimum $(k-1)$-assignment of $Y$
lies in $M$, and by Lemma~\ref{lem:aux}, the minimum $k$-assignments
of $X$ and $Y$ correspond entry by entry.

There are two cases to consider.

In the first case, for some $s\le k-1$ and $t\le k-1$,
the minimum $k$-assignment of $Y$ uses the entries $y_{sk}$ and $y_{kt}$.
These entries correspond to the minimum entry in row $s$ of $X$ outside of
$M$ and the minimum entry in column $t$ of $X$ outside of $M$.
In the second case, the minimum $k$-assignment of $Y$ uses entry $y_{kk}$.
This entry corresponds to the minimum entry in the submatrix of $X$
complementary to $M$.
In both cases,
by Proposition~\ref{pro:prob}, the locations of these minima in $X$
are independent of the random variables making
up the entries of $Y$, and thus independent of the events that the
minimum $(k-1)$-assignment of $Y$ lies in $M$ and the minimum
$k$-assignment of $Y$ uses particular entries outside of $M$.

Thus in the first case, the probability that the
minimum entry in the part of row $s$ outside of $M$ comes from column $j'$
is $\frac{r_sc_{j'}}{\sum_{j=k}^n r_sc_j} = \frac{c_{j'}}{c_k+\ldots + c_n}$,
and the probability that the
minimum entry in the part of column $t$ outside of $M$ comes from row $i'$
is $\frac{r_{i'}c_t}{\sum_{i=k}^m r_ic_t} = \frac{r_{i'}}{r_k+\ldots r_m}$.
Moreover, the locations of these minima within row $s$ and column $t$ are independent
events since the parts of row $s$ and column $t$ outside of $M$ are disjoint.

In the second case,
the probability of the minimum entry in the submatrix of $X$ complementary
to $M$
coming from row $i'$ and column $j'$ is $\frac{r_{i'}c_{j'}}{\sum r_ic_j}$
where the sum in the denominator is over all locations $(i,j)$ in the
submatrix of $X$ complementary to $M$.
Thus, the probability that the minimum entry comes from row $i'$ is
$\frac{r_{i'}}{\sum_{i=k}^m r_i}$ and the probability that the minimum
entry comes from column $j'$ is $\frac{c_{j'}}{\sum_{j=k}^n c_j}$.
$\Box$\medskip

From this lemma we can immediately conclude the following theorem,
which imparts further meaning to the formal $\Pr(r,I)$ and
$\Pr(r,(i_1,\dots,i_k))$ functions used in Section~\ref{sec:equiv}:
\begin{theorem}\label{th:flag}
Suppose that $A=(r_ic_j)$ is a rank 1 rate matrix and
$X$ an exponential random matrix with rate matrix
$A$.  Let $(i_1,\dots,i_k)$ be a sequence of distinct elements
of $[m]$ and $(j_1,\dots,j_k)$ a sequence of 
distinct elements from $[n]$.  Then the probability that $X$
has the associated flag of submatrices is 
$\Pr(r,(i_1,\dots,i_k))\Pr(c,(j_1,\dots,j_k))$.
Furthermore, if $I\subseteq [m]$ and $J\subseteq [n]$ are sets of size 
$k$, then the probability that the minimum $k$-assignment
of $X$ uses the rows indexed by $I$ is $\Pr(r,I)$ and
the probability that it uses the columns indexed by $J$ is
$\Pr(c,J)$, and these events are independent.
\end{theorem}
$\Box$\medskip

We will see that this formula enters in an essential way into
the proof of Theorem~\ref{th:elimit}.

\section{Computational evidence for our conjectures}\label{sec:comp}

By Proposition~\ref{pro:prob}, it is an easy matter to compute
the expected value of the minimum 1-assignment for an arbitrary rate
matrix.

\begin{example}
The expected value of the minimum 1-assignment when the rate matrix is
$A=(a_{ij})$ is
$$\frac{1}{\sum_{ij}a_{ij}}.$$ 
\end{example}

For $k\ge 2$ the computation is more complicated.
In \cite{as} and \cite{cs}
the authors calculate the expected value
of the minimum assignment for a random exponential matrix when the
rates are all 1 and $k$ is small.

The method in \cite{as} applies just as well to the
case of arbitrary rate matrices. 
The essence of their idea is to introduce a slightly more general 
expectation problem in which they choose all the entries of the random
matrix $X$ as before, except that there is a set $Z$ of fixed zeroes
in $X$.  Let us denote the expected value of the minimum assignment
in this case by $E(A,Z)$.  

It is then sometimes possible to establish a recursive calculation
of $E(A,Z)$.  The base of the recursion occurs when there exist
$k$ zeroes in  $Z$, no two in the same row or column.  In
this case we know that the expected value of the minimum assignment
is zero.  For the inductive
part of the calculation we can sometimes express an expected value
$E(A,Z)$ as a constant plus a 
linear combination of $E(A,Z')$ where $Z'$ is obtained from $Z$ 
by adjoining one more position to $Z$.

This arises as follows. Suppose that $X$ is a random exponential
matrix except for a set $Z$ of positions in $X$ where the entries are
fixed zeroes.  Suppose further that we have a set $S$ of positions in
$X$, disjoint from $Z$, such that any minimum $k$-assignment of $X$
meets the set $S$ in exactly $r$ positions. (In other words, every
nonnegative matrix with zero set $Z$ has the property that its minimum
$k$-assignments all meet $S$ in exactly $r$ positions.) Abusing
notation, we also let $S$ denote the matrix which is 1 at the
positions in the set $S$ and zero otherwise.

We will derive the following formula:
\begin{equation}\label{eq:4}
E(A,Z)=\frac{r}{A\cdot S}+\sum_{(i,j) \in S} \frac{a_{ij}}{A \cdot S} E(A,Z \cup \{(i,j)\} ).
\end{equation}
Indeed, the integral for $E(A,Z)$, which involves only the variables
$x_{ij}$ for $(i,j) \notin Z$, is given by 
\[
E(A,Z)=\left(\prod_{(i,j) \notin Z}a_{ij}\right) \int_{X}{\min}_k(X)
e^{-A \cdot X} dX. 
\]
We can derive (\ref{eq:4}) by breaking up this
integral into $|S|$ parts, each corresponding to a position in $S$
containing the minimum entry among all positions in $S$.  For the part
of the integral where $x_{i_0j_0}$ is the minimum entry in $S$, we
make a change of variables with Jacobian 1, as follows.  We express
the $x_{ij}$ in terms of new variables $y_{ij}$ by setting
$x_{ij}=y_{i_0j_0}+y_{ij}$ when $(i,j)\in S-\{(i_0,j_0)\}$
and $x_{ij}=y_{ij}$ otherwise.  $X$ can then be written as
$Y+y_{i_0j_0}S$ where $Y$ is a nonnegative matrix with fixed zeroes
at $Z \cup \{(i_0,j_0)\}$.  From our hypothesis about $S$, we have
\[
{\min}_k(X)={\min}_k(Y)+r y_{i_0j_0}.
\]
(Otherwise, there would be
a non-minimum $k$-assignment of $X$, meeting $S$ in fewer than $r$
positions, that becomes a minimum $k$-assignment of a matrix $X - tS$
for some $t< y_{i_0j_0}$. But, the matrix $X-tS$ still has zero set
$Z$, so our hypothesis on $S$ would be contradicted.)
Thus, this part of the integral becomes
$$\left(\prod_{(i,j) \notin Z }a_{ij}\right) 
\int_{Y,y_{i_0j_0}}( {\min}_k(Y)+r y_{i_0j_0} )
   e^{-A \cdot (Y+ y_{i_0j_0}S)} dy_{i_0j_0} dY.$$
This can be computed as the sum of two integrals in the obvious way.  The
first is $a_{i_0j_0}E(A,Z \cup \{(i_0,j_0)\})/(A \cdot S)$ and the
second is $ ra_{i_0j_0}/(A \cdot S)^2$.  When we sum these expressions
over all $(i_0,j_0) \in S$ we obtain (\ref{eq:4}).

When $k=m=n \le 4$, it is easy to see that, when we are not in
the base case, there always exists a set $S$ of positions in $X$ and 
disjoint from $Z$ such that every minimum $k$-assignment of $X$
meets $S$ in the same number of positions.  

To illustrate the method we now discuss the case $k=m=n=4$.
First note that if any row or column of $X$ has no fixed zeroes,
then we can take that row or column to be the set $S$.
So we can suppose that every row or column has at least
one fixed zero.

Now suppose that there is a $3 \times 3$ submatrix $S$ of $X$ 
that has no fixed zeroes.  Without loss of generality, we
may take this to be the upper left  $3 \times 3$ submatrix, so
the matrix $X$ has the form 
\begin{equation}
X=\left[\begin{array}{rrrr}
   * & * & * & 0 \\
   * & * & * & 0 \\
   * & * & * & 0 \\
   0 & 0 & 0 & . \\
\end{array}\right]
\end{equation}
where $*$ means that entry is positive
and $.$ means that nothing is known about that entry.
Then any minimum 4-assignment must use either
two or three entries from $S$. If it uses three, then it must
use $x_{44}$ and some entry $x_{ij}$, $i,j<4$.
But then we can decrease the value of the
assignment by replacing $x_{ij}$ and $x_{44}$ with 
$x_{i4}$ and $x_{4j}$, both of which are zero.
This contradicts the minimality of the 4-assignment we
started with. Thus, any minimum 4-assignment must use
exactly two entries from $S$.

If every $3 \times 3$ submatrix of $X$ has a fixed zero, and we are
not in the base case, then the Hall marriage theorem implies that
there is a $2\times 3$ or $3\times 2$ submatrix $S$ that has no fixed
zeroes.  Suppose the former, which we can take to be the upper left
$2\times 3$ submatrix of $X$.  Each of the first three columns has at
least one fixed zero.  The fixed zeros in those columns must be in
more than one row, since every $3\times 3$ submatrix has a fixed zero.
Thus, we may assume the matrix $X$ has the form
\begin{equation}
X=\left[\begin{array}{rrrr}
   * & * & * & 0 \\
   * & * & * & 0 \\
   0 & 0 & . & . \\
   . & . & 0 & . \\
\end{array}\right]
\end{equation}
Any minimum $4$-assignment must use one or two entries from $S$. 
Suppose a minimum $4$-assignment uses two entries from
$S$.  It cannot use $x_{11}$, since then there would be a smaller
$4$-assignment consisting of $x_{11}$,$x_{24}$,$x_{32}$,and $x_{43}$.
Similarly it cannot use $x_{12}$, $x_{21}$, or $x_{22}$.  But, it can
only use one of $x_{13}$ and $x_{23}$, so it must use exactly one entry
from $S$. 

Thus, we can always find a suitable set $S$ to continue the recursive
calculation.

We have used this method to compute the expected minimum $k$-assignment
for various small cases.  For the case $k\le m=n\le 3$, this is easily
carried out by Mathematica and confirms our conjecture.

When $k=m=n=4$ and $k=m=n=5$ we were not patient enough to wait for
Mathematica to simplify the complete rational expression, even when the
rate matrix has rank 1.  However, we were able to check that we
obtained the correct answer for many random choices or $r_i$'s and
$c_j$'s.  For this purpose, we used an ordinary C program, but, instead
of using exact rational arithmetic, we carried out our calculations
modulo a large prime.  Even so, the evidence seems to be overwhelming
that our conjecture is correct in these cases.

It is possible, although somewhat more complicated, to compute the
expected contribution of a $(k-1)\times (k-1)$ submatrix to the
expected minimum assignment of a $k \times k$ matrix when $2 \le k \le
4$.  In the cases $k=2$ and $k=3$ we were able to check directly with
Mathematica that Conjecture~\ref{conj:contrib} was valid, which proves
Conjecture~\ref{conj:cs} and Conjecture~\ref{conj:r1} whenever $k\le
3$.

When $k=4$ we obtained computational evidence for the validity of
Conjecture~\ref{conj:contrib}, checking its validity in 
a large number of random cases modulo a prime.  
This provides confirmation of the other 
conjectures when $k=4$ and $m$ and $n$ are arbitrary.

\section{Additional evidence for the main conjecture}\label{sec:extra}

Let $A=(r_ic_j)$ as usual and denote $E_k(A)$ by
$E(k,r_1,\dots,r_m,c_1,\dots,c_n)$, or simply $E(k,r,c)$.
Recall from Section~\ref{sec:equiv} 
the notation $F(k,r,c)$ for the formula in Conjecture~\ref{conj:r1}.
In this section we will show that $E$ and $F$ share several properties.

Let us consider how $E$ behaves if we let a collection of the
$r_i$'s approach 0.  For simplicity we assume that $r_1,r_2,\dots,r_l$
approach 0.  Then the random matrices will have very large entries in
the first $l$ rows.  When $k \le m-l$, there are
assignments which avoid the first $l$ rows, so in the limit
that $r_1,r_2,\dots,r_l \rightarrow 0$, 
the minimum assignment will avoid those rows and become equal to
$ E(k,r_{l+1},\dots,r_m,c).$

Now suppose that $k>m-l$.  Then a $k$-assignment must use 
at least $k-(m-l)=k+l-m$ of the first $l$ rows.  But, in our limiting case,
these rows will be very large, so the minimum $k$-assignment will use
as few as possible, or exactly $k+l-m$ of them.
The contribution of the entries from these rows to the 
minimum $k$-assignment will dominate the minimum $k$-assignment, so
in the minimum $k$-assignment this contribution will be as small as possible.
In particular, in the
limit as $r_1,r_2,\dots,r_l \rightarrow 0$, this part of the minimum
$k$-assignment will be 
$E(k+l-m,r_1,\dots,r_l,c_1,\dots,c_n)$.  By Theorem~\ref{th:flag}
we know that a set $K$ of $k+l-m$ columns will be
used by the part of the assignment in the first $l$ rows
with probability $\Pr(c,K)$.  When this happens the expected contribution
from the remaining rows is
$E(m-l,r_{l+1},\dots,r_m,c'(K))$
where by $c'(K)$ denotes the $c_j$'s corresponding to columns not in $K$.
Thus, we should have the following
\begin{theorem}\label{th:elimit}
When $k \le m-l$,
\begin{equation}\label{eq:elimit1}
\lim_{r_1,\dots,r_l \rightarrow 0}E(k,r,c)
=  E(k,r_{l+1},\dots,r_m,c).
\end{equation}
When $k>m-l>0$, 
\begin{flalign}\label{eq:elimit2}
\lim_{r_1,\dots,r_l \rightarrow 0} (E(k,r,c)&-E(k+l-m,r_1,\dots,r_l,c))\notag  \\
&= \sum_K \Pr(c,K) E(m-l,r_{l+1},\dots,r_m,c'(K)) 
\end{flalign}
where the sum is over $K \subseteq [n]$ such that $|K|=k+l-m$.
\end{theorem}
{\em Proof:\/}\ Let $Z$ be a random exponential $m\times n$ matrix
with all entries of mean 1. Then define $X = Z/A$ to be the term by
term quotient of the random matrix $Z$ by the fixed rate matrix $A$,
where $a_{ij}=r_ic_j$. Then $X$ is a random exponential matrix with
rate matrix $A$. In particular,
\[
E(k,r,c)=E_k(A)= E({\min}_k(Z/A)).
\]

Let $A_u$ denote the $ l\times n$ matrix that is comprised of the first $l$ 
rows of $A$, and let $A_d$ denote $ (m-l)\times n$ matrix consisting of the last
$m-l$ rows of $A$. Furthermore, let $Z_u$ and $Z_d$ denote random
exponential matrices of the same corresponding shapes, again with rate 1.

If $k \leq m-l$, then it is easy to see that
\[
\lim_{r_1,\dots,r_l \rightarrow 0} {\min}_k(Z/A) = {\min}_k(Z_d/A_d)
\]
pointwise almost everywhere (i.e., \emph{almost surely}, as random variables).
Since ${\min}_k(Z/A) \leq {\min}_k(Z_d/A_d)$, and the expectation of
${\min}_k(Z_d/A_d)$ is finite, we can apply the dominated convergence
theorem to show that the limit of the expectation is equal to the
expectation of the limit. Thus,
\[
\lim_{r_1,\dots,r_l \rightarrow 0} E({\min}_k(Z/A)) =
E({\min}_k(Z_d/A_d))
\]
and therefore
\[
\lim_{r_1,\dots,r_l \rightarrow 0} E(k,r,c) =
E(k,r_{l+1},\dots,r_m,c).
\]

Now suppose $k>m-l>0$.  Consider the nonnegative random variable $
R ={\min}_k(Z/A) - {\min}_{k+l-m} (Z_u/A_u) $. (It is nonnegative because
the minimum $k$-assignment of $Z/A$ must use at least $k+l-m$ elements
from $Z_u/A_u$.) 
Furthermore, 
\[
{\min}_k(Z/A) \leq {\min}_{k+l-m} (Z_u/A_u)+ {\max}_{m-l} (Z_d/A_d)
\]
where ${\max}_k X$ denotes the \emph{maximum} $k$-assignment of a
matrix $X$. Then the second summand is a nonnegative random variable of
finite expectation, independent of $r_1,\dots,r_l$, dominating $R$.
Let $\chi_K$ denote the random variable that is 1 or 0 depending upon
whether the minimum $(k+l-m)$-assignment of the matrix $Z_u/A_u$ uses
precisely the columns from the set $K$ or does not.
Then, by Theorem~\ref{th:flag}, we know that $E(\chi_K)= \Pr(c,K)$,
independent of $r$.  Let $C_K(Z_d/A_d)$
denote the submatrix of $Z_d/A_d$ obtained when the columns indexed by
$K$ are removed.  We can then show
\[
\lim_{r_1,\dots,r_l \rightarrow 0}\left( R - \sum_K \chi_K {\min}_{m-l}(C_K(Z_d/A_d))\right) = 
0
\]
with convergence pointwise almost everywhere, where the sum is over $K\subseteq [n]$
such that $|K|=k+l-m$. Furthermore, we can bound the finite sum by a random
variable independent of $r_1,\dots,r_l$. (The random variable $\sum_K
{\min}_{m-l}(C_K(Z_d/A_d))$ will do.) Thus, by the dominated
convergence theorem, we can take the limit of the expectations and
obtain
\begin{eqnarray*}
\lim_{r_1,\dots,r_l \rightarrow 0} E(R) &=& \lim_{r_1,\dots,r_l
\rightarrow 0} \sum_K E(\chi_K {\min}_{m-l}(C_K(Z_d/A_d)))  \\
&=& \lim_{r_1,\dots,r_l \rightarrow 0} \sum_K E(\chi_K) E(
{\min}_{m-l}(C_K(Z_d/A_d))) \\
&=&  \sum_K \Pr(c,K) E(m-l,r_{l+1},\dots,r_m,c'(K))
\end{eqnarray*}
But, $ E(R)=
E(k,r,c)-E(k+l-m,r_1,\dots,r_l,c)$, so we are done.$\Box$\medskip

We stated Theorem~\ref{th:elimit} in terms of limits as the first
$l$ $r$'s approach 0, in order to simplify the notation.  However,
since $E(k,r,c)$ is symmetric in the $r$'s and $c$'s, the
analogous results hold for any set of $l$ $r$'s approaching 0.

From Example~\ref{ex:1}, we have that
\begin{equation}\label{eq:aa0}
E(1,r_1,\ldots,r_m,c_1,\ldots,c_n)=\frac{1}{(\sum_i r_i)(\sum_j c_j)}. 
\end{equation}
Also when $l=1$, taking symmetry into account, Theorem~\ref{th:elimit} 
reduces to the following.

Suppose that $m>1$ and $1 \le i \le m$.  Then, if $k<m$, 
\begin{equation}\label{eq:aa1}
\lim_{r_i \rightarrow 0}E(k,r,c)
=  E(k,r_1,\ldots,\hat r_i,\dots,r_m,c) \end{equation}
while if $k=m$, 
\begin{flalign}\label{eq:aa2}
\lim_{r_i \rightarrow 0} 
\bigg( & E(k,r,c)-\frac{1}{r_i \sum_j c_j } \bigg)
\notag \\
&= \frac
{\sum_j c_j E(k-1,r_1,\ldots,\hat r_i,\ldots,r_k,c_1,\ldots,\hat c_j,\ldots,c_n)}{\sum_j c_j}.
\end{flalign}

\begin{proposition}\label{prop:2}
There is at most one set of functions 
$$G(k,r_1,\ldots,r_m,c_1,\ldots,c_n),$$
each a linear combination of terms
of the form~(\ref{eq:recipsum}),
that satisfy the equations (\ref{eq:aa0}), (\ref{eq:aa1}) and (\ref{eq:aa2})
with $E$ replaced by $G$.
\end{proposition}
{\em Proof:\/}\
Let $H(k,m,n)=H(k,r_1,\ldots,r_m,c_1,\ldots,c_n)$ 
denote the difference between two sets of functions satisfying the
conditions of the proposition.  We will show by induction on $k$ and $m$ that
$H(k,m,n)=0$.

It is clear that $H(1,m,n)=0$ for all $m$ and $n$.   Given values of $k$ and $m$,
we may suppose that $H(k',m',n)=0$  and
$H(k,m',n)=0$ for $k'<k$, $m'<m$, and arbitrary $n$.
Then, by equations (\ref{eq:aa1}) and (\ref{eq:aa2})
and induction, we have, for any $i$,  $\lim_{r_i \rightarrow 0} H(k,m,n)=0$.
Suppose that 
$$H(k,m,n)=\sum_I \frac{h_I}{\sum_{i\in I}r_i},  $$
where $I$ runs over nonempty subsets of $[m]$ and the $h_I$'s are 
rational functions of the $c$'s.  Suppose that $h_I \ne 0$ for
some $I$.  Let $I_0$ be a minimal $I$ such that 
$h_{I} \ne 0$ and let $t \in I_0$.
Since $\lim_{r_t \rightarrow 0} H(m,n,k)$ exists, we must have 
$h_{\{t\}}=0$, so $I_0 \ne \{t\}$.
Also
$$  \lim_{r_t \rightarrow 0} H(k,m,n)=
\sum_I 
 \frac{h_I+h_{I-\{t\}}}{\sum_{i \in I-\{t\}}r_i} 
$$
where the sum is over all $I$ strictly containing $\{t\}$. Since the
terms $\frac{1}{\sum_{i \in I}r_i}$ are linearly independent, we must
have $h_I+h_{I-\{t\}}=0$ for all $I$ strictly containing $\{t\}$.
In the case  $I=I_0$, this contradicts the minimality of $I_0$.
Thus $h_I=0$ for all $I$ and $H(k,m,n)=0$.
$\Box$\medskip

Now we show that the rational functions  $F(k,r_1,\ldots,r_m,c_1,\ldots,c_n)$ 
satisfy the same limit conditions that are proved about 
$E(k,r_1,\ldots,r_m,c_1,\ldots,c_n)$ in Theorem~\ref{th:elimit}.

In particular the functions $F(k,r,c)$ satisfy the conditions of Proposition~\ref{prop:2}.
Thus they are the only possible linear combinations of 
terms of the form~(\ref{eq:recipsum}) that could equal $E(k,r,c)$.

The fact that
the $F$'s satisfy {\em all} the limit conditions proved about $E$ provides
additional evidence for Conjecture~\ref{conj:r1}. 

\begin{theorem}\label{th:flimit}
When $k \le m-l$,
\begin{equation}\label{eq:flimit1}
\lim_{r_1,\dots,r_l \rightarrow 0}F(k,r,c)
=  F(k,r_{l+1},\dots,r_m,c).
\end{equation}
When $k>m-l>0$, 
\begin{flalign}\label{eq:flimit2}
\lim_{r_1,\dots,r_l \rightarrow 0} (F(k,r,c)&-F(k+l-m,r_1,\dots,r_l,c))\notag  \\
&= \sum_K \Pr(c,K) F(m-l,r_{l+1},\dots,r_m,c'(K)) 
\end{flalign}
where the sum is over $K \subseteq [n]$ such that $|K|=k+l-m$.
\end{theorem}
{\em Proof:\/}\ 
We use the alternate form of $F(k,r,c)$ given by (\ref{eq:formula1}):
\begin{equation*}
 F(k,r,c)=
\sum_{I,J}
   \binom{k-1-m-n}{k - 1 - |I| - |J|}
      \frac{1}{( \sum_{i \notin I}r_i)  ( \sum_{j \notin J} c_j )},
\end{equation*}
where the binomial coefficient enforces the condition
$|I|+|J|<k$.  

Now, on the left side of 
(\ref{eq:flimit2}), before passing to the limit,
the first term is
$$ \sum_{I,J} 
   \binom{k-1-m-n}{k - 1 - |I| - |J|}
      \frac{1}{( \sum_{i \notin I}r_i)  ( \sum_{j \notin J} c_j )} $$
where the sum is over subsets $I\subseteq [m]$ and $J\subseteq [n]$
and the second term is
$$ \sum_{I,J} 
   \binom{k+l-m-1-l-n}{k+l-m - 1 - |I| - |J|}
      \frac{1}{( \sum_{i \notin I,i\le l}r_i) ( \sum_{j \notin J} c_j )}$$
where the sum is over subsets $I\subseteq [l]$ and $J\subseteq [n]$.
By substituting $I\cup \{l+1,\dots,m\}$ for $I$ in the second term we get
$$ \sum_{I,J}
   \binom{k-1-m-n}{k - 1 - |I| - |J|}
      \frac{1}{( \sum_{i \notin I}r_i) ( \sum_{j \notin J} c_j )}$$
where the sum is over subsets $I\subseteq [m]$ and $J\subseteq [n]$
such that $\{l+1,\dots,m\}\subseteq I$.
Thus, before taking the limit, the left side of (\ref{eq:flimit2}) equals
\begin{equation}\label{eq:ex1}
 \sum_{I,J} 
   \binom{k-1-m-n}{k - 1 - |I| - |J|}
      \frac{1}{( \sum_{i \notin I}r_i) ( \sum_{j \notin J} c_j )}.
\end{equation}
where the sum is over 
all subsets  $J\subseteq [n]$
and those subsets $I\subseteq [m]$ 
that do not contain $\{l+1,\dots,m\}$. 

The last condition on $I$ implies that $\sum_{i\notin I}r_i$ is
nonzero if we set $r_1,\dots,r_l=0$.  Thus,
we can obtain the limit on the left of 
(\ref{eq:flimit2}) simply by replacing $r_1,\dots,r_l$ by zero in
(\ref{eq:ex1}).  When we do this, the effect is that we combine terms with
$I$ having a fixed intersection with $\{l+1,\dots,m\}$.  Note
that this intersection is never all of $\{l+1,\dots,m\}$.

Suppose then that $I$ is a set strictly contained in $\{l+1,\dots,m\}$.
For each $i$ 
there are $\binom{l}{i}$ ways of extending $I$ to a $(|I|+i)$-element
subset of $[m]$ whose intersection with $\{l+1,\dots,m\}$ is $I$.
Thus, after taking the limit on the left of (\ref{eq:flimit2}), we obtain 
\begin{eqnarray*}
&& \sum_{I,J}  \left( \sum_{i=0}^l \binom{l}{i}
      \binom{k-1-m-n}{k-1-|I|-|J|-i} \right) 
  \frac{1}{( \sum_{i \notin I, i > l}r_i)  ( \sum_{j \notin J} c_j )} \\
&=& \sum_{I,J} 
      \binom{k+l-1-m-n}{k-1-|I|-|J|}
  \frac{1}{( \sum_{i \notin I, i > l}r_i) ( \sum_{j \notin J} c_j )} \\
\end{eqnarray*}
where the sum is over proper subsets $I\subsetneq \{l+1,\dots,m\}$ and
all subsets $J\subseteq [n]$.
Note that in the case that $k \le m-l$, this expression is precisely 
$F(k,r_{l+1},\dots,r_m,c)$ so we have proved (\ref{eq:flimit1}).

We continue with the proof of (\ref{eq:flimit2}).
We obtain a slightly more convenient expression if we replace $I$ by 
$I\cup [l]$ in the preceding expression.  Then the left side becomes
\begin{equation}
\sum_{I,J} 
      \binom{k+l-1-m-n}{k+l-1-|I|-|J|}
  \frac{1}{( \sum_{i \notin I}r_i)  ( \sum_{j \notin J} c_j )} 
\end{equation}
where the sum is over sets $I$ strictly contained in 
$[m]$ and containing $[l]$ and
$J\subseteq [n]$.

Now we turn to the right side of (\ref{eq:flimit2}). 

The expression (\ref{eq:formula1}) gives $F(m-l,r_{l+1},\dots,r_m,c'(K))$
as a sum over certain subsets $I$ of $\{l+1,\dots,m\}$ and $J$ of
$c'(K)$.  But this expression is simpler if we replace
$I$ by $I\cup [l]$ and $J$ by $J\cup K$.
Then the right side can be written
$$ \sum_{I \subseteq [m],K \subseteq J \subseteq  [n]}
   \Pr(c,K)\binom{k+l-1-m-n}{k+l-1-|I|-|J|}
  \frac{1}{( \sum_{i \notin I}r_i) 
           ( \sum_{j \notin J} c_j )}$$
where the sum is over proper subsets $I\subsetneq [m]$ containing
$[l]$, subsets $K\subseteq J\subseteq [n]$
such that $|K|=k+l-m$.

Now we have shown that both the left and right sides of 
(\ref{eq:flimit2}) 
are linear combinations of the same reciprocal sums
$1/( \sum_{i \notin I}r_i)$, so to prove (\ref{eq:flimit2})  it
will suffice to prove that the coefficients of the same reciprocal
sums are equal on both sides.

For a given $I\subsetneq [m]$, the coefficient on the left and right depend only
on the cardinality of $I$.  We introduce the 
abbreviations $H=k+l-1-|I|$ and 
$L=k+l-m$. 
Then $0 < L \leq H <k $.  After using these 
abbreviations and equating coefficients we are reduced to
proving 
\begin{equation}\label{eq:basic}
\sum_{J\subseteq [n]}
\binom{L-n-1}{H-|J|}
  \frac{1}{\sum_{j \notin J} c_j} = \sum_{K \subseteq J\subseteq [n]}
   \Pr(c,K) \binom{L-n-1}{H-|J|}
  \frac{1}{\sum_{j \notin J} c_j }
\end{equation}
where in the sum on the right the subset $K$ must have cardinality
$L$.

Now use the expression (\ref{eq:incex}) for $\Pr(c,K)$ to rewrite the 
right side of (\ref{eq:basic}) as
\begin{eqnarray*}
& &\sum_{A\subseteq K\subseteq J \subseteq [n]} \left(
  (-1)^{L-|A|}\frac{\sum_{j\notin K} c_j}{\sum_{j\notin A} c_j} \right)
 \binom{L-n-1}{H-|J|}
  \frac{1}{\sum_{j \notin J} c_j }\\
\end{eqnarray*}
where we still require that $|K|=L$.
We sum this first over $K$.  In the term
$\sum_{j\notin K} c_j$, the number of $K$'s for which a given $c_j$ occurs
depends only on whether $j$ belongs to $J$.  Thus,
we can rewrite the right side of (\ref{eq:basic}) as
\begin{eqnarray*}
 & &\sum_{A\subseteq J} (-1)^{L-|A|}
\left(\frac{
\binom{|J|-|A|}{L-|A|}\sum_{j\notin J}c_j + 
\binom{|J|-|A|-1}{L-|A|}\sum_{j\in J-A}c_j
}{\left(\sum_{j\notin A}c_j\right) \left(\sum_{j\notin J}c_j \right)} \right)
\binom{L-n-1}{H-|J|}.\\
\end{eqnarray*}
The numerator of the fraction can be rewritten as
\begin{eqnarray*}
&& 
\binom{|J|-|A|}{L-|A|}\sum_{j\notin J}c_j + 
\binom{|J|-|A|-1}{L-|A|}\sum_{j\notin A}c_j -
\binom{|J|-|A|-1}{L-|A|}\sum_{j\notin J}c_j \\
&=&\binom{|J|-|A|-1}{L-|A|-1}\sum_{j\notin J}c_j + 
\binom{|J|-|A|-1}{L-|A|}\sum_{j\notin A}c_j  \\
\end{eqnarray*}
so the right side of (\ref{eq:basic}) can be rewritten as
a sum of two terms:
\begin{equation}\label{eq:term1}
\sum_A (-1)^{L-|A|}\left(
     \sum_{J\supseteq A} \binom{|J|-|A|-1}{L-|A|-1}\binom{L-n-1}{H-|J|}
     \right)
     \frac{1}{\sum_{j\notin A}c_j}
\end{equation}
and
\begin{equation}\label{eq:term2}
\sum_J \left(\sum_{A\subseteq J} (-1)^{L-|A|}\binom{|J|-|A|-1}{L-|A|}\right)
\binom{L-n-1}{H-|J|}\frac{1}{\sum_{j\notin J}c_j}.
\end{equation}

Now, comparing with the left side of (\ref{eq:basic}),
it suffices to show that the inner sum
in (\ref{eq:term1}) equals $(-1)^{L-|A|}\binom{L-n-1}{H-|A|}$ when
$|A|<L$ and 0 otherwise, and the inner
sum in (\ref{eq:term2}) equals 1.

It is easy to see that the inner sum of (\ref{eq:term1}) equals 0
when $|A|=L$, since the first binomial coefficient has a negative
lower term in that case.  When $|A|<L$, because $|J|\ge L$,
the inner sum is over $J$ strictly larger than $A$.
Collecting terms according to the cardinality $j$ of $J$, we obtain 

\begin{flalign*}
 &\quad\quad  \sum_{j={|A|+1}}^{H} \binom{n-|A|}{j-|A|}\binom{j-|A|-1}{L-|A|-1}
                        \binom{L-n-1}{H-j}\\
 &= \sum_{i=1}^{H-|A|} \binom{n-|A|}{i}\binom{i-1}{L-|A|-1}
                        \binom{L-n-1}{H-|A|-i}\\
 &= -\binom{n-|A|}{0}\binom{0-1}{L-|A|-1}\binom{L-n-1}{H-|A|-0} \\ &
\mbox{}\quad \quad +
       \sum_{i=0}^{H-|A|} \binom{n-|A|}{i}\binom{i-1}{L-|A|-1}\binom{L-n-1}{H-|A|-i}\\
 &= -\binom{-1}{L-|A|-1} \binom{L-n-1}{H-|A|} + \binom{n-L+H-|A|}{H-|A|} \binom{-1}{L-1-H}\\
 &= (-1)^{L-|A|}\binom{L-n-1}{H-|A|}\\
\end{flalign*}
where the third equality holds by substituting into the identity
(\cite{r}, p.16)
$$\binom{m}{p}\binom{n}{q}= \sum_{i=0}^p \binom{n+i}{p+q}\binom{m-n+q}{i}\binom{n-m+p}{p-i}$$
and the fourth equality holds because $L-1-H<0.$ 

In the inner sum of (\ref{eq:term2}) we can collect terms according to
the cardinality $a$ of $A$.  Recalling that $A$ must be contained in $J$,
we obtain
\begin{eqnarray*}
\sum_a (-1)^{L-a}\binom{|J|}{a}\binom{|J|-a-1}{L-a}= \sum_a \binom{|J|}{a}\binom{L-|J|}{L-a} 
=  \binom{L}{L}=1.
\end{eqnarray*}

This proves Theorem~\ref{th:flimit}.  $\Box$\medskip

Finally we prove that $E$ and $F$ have another property in common. 

\begin{theorem}\label{thm:mono}
Both $E(k,r,c)$ and $F(k,r,c)$ are monotonically decreasing functions
of $r$ and $c$. In particular, if $r_1<r_1'$, then 
\[
E(k,r_1,r_2,\dots,r_m,c) > E(k,r_1',r_2,\dots,r_m,c)
\]
and
\[
F(k,r_1,r_2,\dots,r_m,c) > F(k,r_1',r_2,\dots,r_m,c).
\]
Furthermore, $F$ is differentiable to any degree $\ell\geq 0$ in each $r_i$ and 
$c_j$, and $ (-1)^\ell\frac{\partial^\ell F}{\partial r_i^\ell}>0$. 
\end{theorem}
{\em Proof:\/}\ Without loss of generality, because of symmetry, we
can restrict ourselves to considering the behavior of $E$ and $F$ as
functions of $r_1$. 

Recall from the proof of Theorem~\ref{th:elimit} that $E(k,r,c)$ is
the expectation of the random variable ${\min}_k (Z/A)$, where $Z$ is
an $m\times n$-matrix-valued random variable with exponentially
distributed independent entries of mean 1, where $A=(r_ic_j)$ is the
rank 1 rate matrix, and where $Z/A$ denotes the element by element
quotient.  Suppose $r_1<r_1'$ and $r_i=r_i'$ for
$i=2,\dots,m$. Define $A'=(r_i'c_j)$. Then $A \leq A'$ term by term,
so $Z/A \geq Z/A'$, and $ {\min}_k (Z/A) \geq {\min}_k (Z/A')$. Hence,
\[
E(k,r,c)= E({\min}_k (Z/A)) \geq E({\min}_k (Z/A')) = E(k,r',c).
\]
Since there is a nonzero probability that the minimum
$k$-assignment of $Z/A$ uses the first row, the inequality is actually
strict.

We start with the following formula for $F$, from (\ref{eq:formula2})
and (\ref{eq:inex}):
\begin{equation*}
F(k,r,c)= \sum_{ |I'|+|J'| < k}
\left( \sum_{I \subseteq I'} (-1)^{|I'|-|I|}\frac{1}{\sum_{i \notin I}r_i}\right)
 \left( \sum_{J \subseteq J'} (-1)^{|J'|-|J|}\frac{1}{\sum_{j \notin J}c_j}\right)
\end{equation*}
For $\ell\geq 1$, $ I'\subsetneq [m]$, and $ J'\subsetneq [n]$, we define the functions
\begin{eqnarray*}
f(\ell,r,I') &=& \sum_{I \subseteq I'} (-1)^{|I'|-|I|}\frac{1}{(\sum_{i
\notin I}r_i)^\ell} \\
g(\ell,c,J') &=& \sum_{J \subseteq J'} (-1)^{|J'|-|J|}\frac{1}{(\sum_{j
\notin J}c_j)^\ell},
\end{eqnarray*}
so that
\begin{equation}\label{eq:Ffg}
F(k,r,c)= \sum_{ |I'|+|J'| < k} f(1,r,I')g(1,c,J').
\end{equation}

First we prove that $f(\ell,r,I')>0$. We define the partial sum
$R=\sum_{i\notin I'}^m r_i$.  Then
\begin{eqnarray*}
0&<& \int_0^\infty t^{\ell-1} e^{-Rt} \left(\prod_{i\in I'} (1-e^{-r_it})\right) dt \\
&=& \int_0^\infty t^{\ell-1} \sum_{I\subseteq I'} (-1)^{|I|}\exp(t(-R-\sum_{i\in I} 
r_i)) dt \\
&=& \int_0^\infty t^{\ell-1} \sum_{I\subseteq I'} (-1)^{|I'|-|I|}\exp(t(-\sum_{i\notin I} 
r_i)) dt \\
&=& \sum_{I\subseteq I'} (-1)^{|I'|-|I|} \int_0^\infty t^{\ell-1} \exp(t(-\sum_{i\notin I} 
r_i)) dt \\
&=& (\ell-1)! \sum_{I\subseteq I'} (-1)^{|I'|-|I|} \frac1{(\sum_{i\notin I} 
r_i)^\ell} \\
&=& (\ell-1)! f(\ell,r,I'),
\end{eqnarray*}
so $f(\ell,r,I')>0$ as claimed.
The proof that $g(\ell,c,J')>0$ is essentially the same.
 
To finish the proof of the theorem it suffices to show that for $\ell\geq1$,
\begin{equation}\label{eq:diffF}
(-\frac{\partial}{\partial r_1})^\ell F(k,r,c)= \sum_{ |I'|+|J'| =
k-1,\ 1\notin I'} \ell! f(\ell+1,r,I')g(1,c,J'),
\end{equation}
because all terms on the right hand side are positive.
First, rewrite (\ref{eq:Ffg}) as follows
\begin{eqnarray*}
F(k,r,c) &=& \sum_{ |I'|+|J'| < k-1, 1\notin I'}
(f(1,r,I')+f(1,r,I'\cup \{1\}))g(1,c,J') \\
&+& \sum_{ |I'|+|J'| = k-1, 1\notin I'}
f(1,r,I')g(1,c,J').
\end{eqnarray*}
The functions $f(1,r,I')+f(1,r,I'\cup \{1\})$ and $g(1,c,J')$
are independent of $r_1$, so those partial derivatives vanish.
Meanwhile, if $i\notin I'$,
then all the denominators in $f(1,r,I')$ involve $r_1$, so we have
$$ (-\frac{\partial}{\partial r_1})^\ell f(1,r,I') = \ell! f(\ell+1,r,I').$$
This proves (\ref{eq:diffF}) and Theorem~\ref{thm:mono}.  $\Box$\medskip

The functions $E(k,r,c)$ and $F(k,r,c)$ share many properties. 
Both functions are rational, homogeneous of degree $-1$, and symmetric
 in the $r_i$'s and the
$c_j$'s. Both rational functions have denominators that factor into
linear factors which are either sums of $r_i$'s or sums of
$c_j$'s. Both are positive, and both are monotonically decreasing
in each $r_i$ and $c_j$.
Finally, $E$ and $F$ share various limit properties
with $r_i$ or $c_j$ tending to 0.
We could also consider limits as $m$ or $n$ go to $\infty$.
However, it is conceivable that the properties we have already found
are sufficient to guarantee that such a function is unique. In
any case, one plan to prove $E=F$ would be to extend the list of
common properties until equality is forced.

\end{document}